# On the Jacobian Question

Dhananjay P. Mehendale
Sir Parashurambhau College, Tilak Road, Pune 411030
India

## Abstract

The direct or algorithmic approach for the Jacobian problem, consisting of the direct construction of the inverse polynomials is proposed. The so called principle and derived Jacobi conditions are proposed and discussed. The algorithmic approach is shown to be extendable to higher dimensions by proceeding on exactly identical lines. As per the important result due to Bass, Connell, and Wright [3] it is enough to show the validity of the Jacobian conjecture for cubic polynomials of special type (BCW form) in two, three, …., n variables. Firstly, the method of proof for the cases of two and three variables is discussed at length. It is then indicated that the extension to the several variables case follows automatically by just following the same steps and there is no hindrance as there is essentially no change in the basic situation and the same line of thought used for the case of two and three variables remains applicable. Thus, we show that the problem can be solved completely using the important reduction of the problem to the case of special cubic degree polynomials [3]. We have shown how to obtain inverse polynomials. We fully obtain them for two variables case and almost obtain them for three variables case.

**1. Introduction:** Given $n$ polynomials $u = (u_1, u_2, \cdots, u_n)$ in $n$ variables $x = (x_1, x_2, \cdots, x_n)$ and their Jacobian $J_x(u) = \det\left(\frac{\partial u_i}{\partial x_j}\right)$, $i, j = 1,2,\cdots,n$ is a nonzero constant in the ground field $k$ of characteristic zero. The problem called the Jacobian problem [1] or the Jacobian conjecture [3], [4] is to show that $k[x_1, x_2, \cdots, x_n] = k[u_1, u_2, \cdots, u_n]$.

In other words, what we call the direct or algorithmic approach, one needs to show that one can construct $n$ polynomials $x = (x_1, x_2, \cdots, x_n)$ in $n$ variables $u = (u_1, u_2, \cdots, u_n)$.

This problem was first proposed for the polynomials in two variables with integral coefficients by O. H. Keller [2].

The Jacobian conjecture has negative answer when the ground field $k$ has a positive characteristic [1].

We begin with few definitions and a simple observation. The observation, though not required to be used in the later course, tells about a property that should be satisfied by the degrees of the polynomials.

**Definition 1.1:** The $x_i$ degree of $u_j$, $\deg_{(x_i)}(u_j)$, is the highest power of $x_i$ appearing in the polynomial $u_j$ with a nonzero coefficient.

Let the monomials $v_l$ in $u_j$ be of type $\alpha(x_1^{s_1} x_2^{s_2} \cdots x_n^{s_n})$.

**Definition 1.2:** The total degree of a monomial of type $v_l$, $\deg(v_l)$, is equal to the sum of all the indices, namely, $\deg(v_l) = \sum_{j=0}^{n} s_j$ .

**Definition 1.3:** The total degree of a polynomial $u_j$ is the maximum of the degrees of monomials present in that polynomial.

**Observation:** Note that when $x_i$ is expressed as inverse function in the variables $u_1, u_2, \cdots, u_n$, then the basic necessity is that one should get back the equation $x_i = x_i$ when all the variables $u_j$ are back substituted in terms of their original form (i.e. as polynomials in $x_1, x_2, \cdots, x_n$). For this to happen all the higher order terms should disappear. In order to get this effect the following conditions must be fulfilled automatically: Thus, if $\deg_{(x_i)}(u_j) = m_j$ for all $j$ then

$$\deg_{(u_j)}(x_i) = p \times \frac{l.c.m.(m_1, m_2, \cdots, m_n)}{\deg_{(x_i)}(u_j)}$$, where $p = 1$, or some positive integer.

**2. A Direct or Algorithmic Approach:** Our approach to the Jacobian problem is constructive, i.e. we propose a systematic procedure to construct the above mentioned $n$ polynomials $x = (x_1, x_2, \cdots, x_n)$ in $n$ variables $u = (u_1, u_2, \cdots, u_n)$. We begin with

**2.1 The Two Variables Case:** Given

$$f = \sum_{m=0}^{k_1} \sum_{n=0}^{k_2} a_{mn} x^m y^n \qquad (2.1.1)$$

and

$$g = \sum_{m=0}^{l_1} \sum_{n=0}^{l_2} b_{mn} x^m y^n \qquad (2.1.2)$$

$a_{mn}, b_{mn} \in k$, and the Jacobian

$$J_{(x,y)}(f,g) = \det\begin{pmatrix} f_x & f_y \\ g_x & g_y \end{pmatrix} = 1 \neq 0 \qquad (2.1.3)$$

We need to show that $x, y$ (which we can express as power series in $f, g$ using the inverse function theorem) are actually polynomials, i.e.

$$x = \sum_{m=0}^{u_1} \sum_{n=0}^{u_2} c_{mn} f^m g^n \qquad (2.1.4)$$

and

$$y = \sum_{m=0}^{v_1} \sum_{n=0}^{v_2} d_{mn} f^m g^n \qquad (2.1.5)$$

$c_{mn}, d_{mn} \in k$.

**2.2 The Local Homeomorphism:** Let $F(x,y) = (f,g)$. Let $f_x, f_y, g_x, g_y$ exist and are continuous, therefore,

$$DF(x,y)=\begin{pmatrix} f_x & f_y \\ g_x & g_y \end{pmatrix}$$ exists and since

$\det(DF(x,y) = J_{(x,y)}(f,g) = 1 \neq 0$, therefore, by the inverse function theorem, $F^{-1}$ exists and moreover

$$DF^{-1}(f,g) = \left(DF(x,y)\right)^{-1},$$

i.e.

$$\begin{pmatrix} x_f & x_g \\ y_f & y_g \end{pmatrix} = \begin{pmatrix} g_y & -f_y \\ -g_x & f_x \end{pmatrix} \qquad (2.2.1)$$

The (local) inverse

$$F^{-1}(f,g) = (x,y) = (\sum_{m=0} \sum_{n=0} c_{mn} f^m g^n, \sum_{m=0} \sum_{n=0} d_{mn} f^m g^n),$$

is made up of power series in general and the Jacobian conjecture demands that they are in fact polynomials (whose coefficients can be uniquely determined in terms of the coefficients of the originally given polynomials $f, g$).

**2.3 Monomials in Variables as a Vector Space Basis:** Suppose we are given the following system of equations:

$$f^* = a_{10}x + a_{01}y \qquad (2.3.1)$$

$$g^* = b_{10}x + b_{01}y \qquad (2.3.2)$$

and $\det\begin{pmatrix} a_{10} & a_{01} \\ b_{10} & b_{01} \end{pmatrix} = 1 \neq 0$. So, clearly we have the following inverse relations

$$x = b_{01}f^* - a_{01}g^* \qquad (2.3.3)$$

$$y = -b_{10}f^* + a_{10}g^* \qquad (2.3.4)$$

**Theorem 2.3.1:** If the system of equations represented by the equations (2.3.1) and (2.3.2) given above is an invertible system with inverse functions given by linear polynomials in (2.3.3) and (2.3.4) then the system of equations formed by constructing the homogeneous blocks of equations of degree $n$, obtained from equations (2.3.1) and (2.3.2) as well as equations (2.3.3) and (2.3.4), is also invertible for any $n$.

**Proof:** Equations (2.3.1) and (2.3.2) can be expressed as the following matrix equation:

$$\begin{pmatrix} f^* \\ g^* \end{pmatrix} = \begin{pmatrix} a_{10} & a_{01} \\ b_{10} & b_{01} \end{pmatrix}\begin{pmatrix} x \\ y \end{pmatrix}$$

and their invertible nature implies that we can write

$$\begin{pmatrix} x \\ y \end{pmatrix} = \begin{pmatrix} b_{01} & -a_{01} \\ -b_{10} & a_{10} \end{pmatrix}\begin{pmatrix} f^* \\ g^* \end{pmatrix}$$

Let $U_2 = \begin{pmatrix} a_{10} & a_{01} \\ b_{10} & b_{01} \end{pmatrix}$ and $V_2 = \begin{pmatrix} b_{01} & -a_{01} \\ -b_{10} & a_{10} \end{pmatrix}$

Since $\{x, y\}$ is an independent set of linear monomials in variables (subspace basis) we have $I - V_2U_2 = 0$, or $V_2 = (U_2)^{-1}$.

Now, by squaring both sides of equations (2.3.1) and (2.3.2) and taking the product of these equations one can construct the homogeneous blocks of equations of degree $2$ as follows which essentially will lead to the construction of the following matrix equation, namely,

$$\begin{pmatrix} (f^*)^2 \\ (f^*g^*) \\ (g^*)^2 \end{pmatrix} = \begin{pmatrix} a_{10}^2 & 2a_{10}a_{01} & a_{01}^2 \\ a_{10}b_{10} & (a_{10}b_{01} + a_{01}b_{10}) & a_{01}b_{01} \\ b_{10}^2 & 2b_{10}b_{01} & b_{01}^2 \end{pmatrix} \begin{pmatrix} x^2 \\ xy \\ y^2 \end{pmatrix}$$

Similarly, by squaring equations (2.3.3) and (2.3.4) and taking the product of these equations one can construct the homogeneous blocks of equations of degree $2$ as follows which essentially will lead to the construction of the following matrix equation, namely,

$$\begin{pmatrix} x^2 \\ xy \\ y^2 \end{pmatrix} = \begin{pmatrix} b_{01}^2 & -2a_{01}b_{01} & a_{10}^2 \\ -b_{10}b_{01} & (a_{10}b_{01} + a_{01}b_{10}) & -a_{01}a_{10} \\ b_{10}^2 & -2a_{10}b_{10} & a_{10}^2 \end{pmatrix} \begin{pmatrix} (f^*)^2 \\ (f^*g^*) \\ (g^*)^2 \end{pmatrix}$$

Thus, let

$$U_3 = \begin{pmatrix} a_{10}^2 & 2a_{10}a_{01} & a_{01}^2 \\ a_{10}b_{10} & (a_{10}b_{01} + a_{01}b_{10}) & a_{01}b_{01} \\ b_{10}^2 & 2b_{10}b_{01} & b_{01}^2 \end{pmatrix}, \text{ and}$$

$$V_3 = \begin{pmatrix} b_{01}^2 & -2a_{01}b_{01} & a_{01}^2 \\ -b_{10}b_{01} & (a_{10}b_{01} + a_{01}b_{10}) & -a_{01}a_{10} \\ b_{10}^2 & -2a_{10}b_{10} & a_{10}^2 \end{pmatrix}$$

Since monomials $\{x^2, xy, y^2\}$ is in independent set of quadratic monomials in variables (subspace basis) therefore, $I - V_3U_3 = 0$, or $V_3 = (U_3)^{-1}$.

By taking all possible products leading to formation of homogeneous block of equations of total degree $n$ by making use of the equations (2.3.1) and (2.3.2), we can form the following matrix equation:

$$\begin{pmatrix} (f^*)^n \\ (f^*)^{n-1}g^* \\ (f^*)^{n-2}(g^*)^2 \\ \vdots \\ (f^*)^2(g^*)^{n-2} \\ f^*(g^*)^{n-1} \\ (g^*)^n \end{pmatrix} = U_n \begin{pmatrix} x^n \\ x^{n-1}y \\ x^{n-2}y^2 \\ \vdots \\ x^2y^{n-2} \\ xy^{n-1} \\ y^n \end{pmatrix}$$

Similarly, using equations (2.3.3) and (2.3.4) and performing the same operations as above we can form the following matrix equation:

$$\begin{pmatrix} x^n \\ x^{n-1}y \\ x^{n-2}y^2 \\ \vdots \\ x^2y^{n-2} \\ xy^{n-1} \\ y^n \end{pmatrix} = V_n \begin{pmatrix} (f^*)^n \\ (f^*)^{n-1}g^* \\ (f^*)^{n-2}(g^*)^2 \\ \vdots \\ (f^*)^2(g^*)^{n-2} \\ f^*(g^*)^{n-1} \\ (g^*)^n \end{pmatrix}$$

Since monomials $\{x^n, x^{n-1}y, \cdots, xy^{n-1}, y^n\}$ is in independent set of n-th degree monomials in variables (subspace basis) therefore, we have $I - V_nU_n = 0$, or $V_n = (U_n)^{-1}$.

□

**Remark 2.3.1:** To solve the Jacobian conjecture for two variables, as already mentioned above, our objective is to express the coefficients of the power series representing inverse functions in terms of the coefficients of the originally given polynomials and to show further that only finitely many of these coefficients of the power series are actually nonzero, i.e. to show that these power series are actually polynomials. In the process of achieving this task we come across similar matrix equations like those described in the proof of the above theorem, namely, the matrix equations of the type $WF = X$, where $W = (U_n)^T$, the transpose of $U_n$, and $F$ is unknown vector of coefficients of power series for $x, y$ to be determined, through the inverse relation $F = (V_n)^T X$, where $X$ is known vector made up in terms of the coefficients of $f, g$.

**2.4 The Jacobi Conditions:** The polynomial pairs $(f, g)$ with coefficients in the field $k$ of characteristic zero (e.g. $k = R$, the field of real numbers, or $k = C$, the field of complex numbers) for which the Jacobian

$$J_{(x,y)}(f,g) = \det\begin{pmatrix} f_x & f_y \\ g_x & g_y \end{pmatrix} = 1 \neq 0$$

is described by saying that the polynomial pairs $(f, g)$ satisfy the **Jacobi Condition**. The fulfillment of the Jacobi condition automatically implies the fulfillment of many other conditions. In order to distinguish them from the above mentioned Jacobi condition let us call the above mentioned Jacobi condition the **Principle Jacobi Condition** and call the automatically fulfilled many other conditions which are implied by this principle Jacobi condition the **Derived Jacobi Conditions**.
Using equations (2.1.1) and (2.1.2) in equation (2.1.3)
(1) We get the constant term as,

$a_{10}b_{01} - a_{01}b_{10} = 1$ (the Principle Jacobi condition) (2.4.1)

As the polynomial pair $(f, g)$ is satisfying the principle Jacobi condition therefore the fulfillment of the following derived Jacobi conditions is automatically implied by this principle Jacobi condition.
(2) Collecting the coefficients of $x$, we (automatically) have

$$(2a_{20}b_{01} + a_{10}b_{11}) - (2b_{20}a_{01} + b_{10}a_{11}) = 0 \qquad (2.4.2)$$

(3) Collecting the coefficients of $y$,

$$(a_{11}b_{01} + 2a_{10}b_{02}) - (b_{11}a_{01} + 2b_{10}a_{02}) = 0 \qquad (2.4.3)$$

(4) Collecting the coefficients of $x^2$,

$$(3a_{30}b_{01} + 2a_{20}b_{11} + a_{10}b_{21}) - (3b_{30}a_{01} + 2b_{20}a_{11} + b_{10}a_{21}) = 0 \quad (2.4.4)$$

(5) Collecting the coefficients of $xy$,

$$(2a_{21}b_{01} + a_{11}b_{11} + 4a_{20}b_{02} + 2a_{10}b_{12})$$
$$-(2b_{21}a_{01} + b_{11}a_{11} + 4b_{20}a_{02} + 2b_{10}a_{12}) = 0 \quad (2.4.5)$$

(6) Collecting the coefficients of $y^2$,

$$(a_{12}b_{01} + 2a_{11}b_{02} + 3a_{10}b_{03})$$
$$-(b_{12}a_{01} + 2b_{11}a_{02} + 3b_{10}a_{03}) = 0 \quad (2.4.6)$$

$$\vdots$$

Thus, continuing on the similar lines we can obtain the derived Jacobi conditions by collecting the coefficients for each monomial $x^{k-r}y^r$ in the homogeneous block $\{x^k, x^{k-1}y, x^{k-2}y^2, \cdots, y^k\}$ of monomials of degree $k$ in the Jacobian determinant.

**2.5 Jacobi Conditions for Inverse Function Pair:** From equation (2.2.1) we have

$$\det\begin{pmatrix} x_f & x_g \\ y_f & y_g \end{pmatrix} = 1 \neq 0$$

This implies the principle as well as derived Jacobi conditions for the inverse function pair $(x, y)$ as follows. These conditions can be obtained by replacing

(i) all $a^s$ by $c^s$

(ii) all $b^s$ by $d^s$, and

(iii) all $x^r y^s$ by $f^r g^s$

(1) We get the constant term as,

$c_{10}d_{01} - c_{01}d_{10} = 1$ (the Principle Jacobi condition) (2.5.1)

(2) Collecting the coefficients of $f$ ,

$$(2c_{20}d_{01} + c_{10}d_{11}) - (2d_{20}c_{01} + d_{10}d_{11}) = 0 \qquad (2.4.2)$$

(3) Collecting the coefficients of $g$ ,

$$(c_{11}d_{01} + 2c_{10}d_{02}) - (d_{11}c_{01} + 2d_{10}c_{02}) = 0 \qquad (2.4.3)$$

(4) Collecting the coefficients of $f^2$,

$$(3c_{30}d_{01} + 2c_{20}d_{11} + c_{10}d_{21}) - (3d_{30}c_{01} + 2d_{20}c_{11} + d_{10}c_{21}) = 0$$

––– (2.4.4)

(5) Collecting the coefficients of $fg$ ,

$$(2c_{21}d_{01} + c_{11}d_{11} + 4c_{20}d_{02} + 2c_{10}d_{12})$$
$$-(2d_{21}c_{01} + d_{11}c_{11} + 4d_{20}c_{02} + 2d_{10}c_{12}) = 0$$

––– (2.4.5)

(6) Collecting the coefficients of $g^2$,

$$(c_{12}d_{01} + 2c_{11}d_{02} + 3c_{10}d_{03})$$
$$-(d_{12}c_{01} + 2d_{11}c_{02} + 3d_{10}c_{03}) = 0$$

––– (2.4.5)

$$\vdots$$

**2.6 Jacobi Conditions with Mixed Entries:** Substituting further the expressions for $\{f, g, f^2, fg, g^2 \cdots\}$ in terms of the variables $\{x, y\}$ we get the new derived Jacobi conditions involving $\{a^s, b^s, c^s, d^s\}$ together, obtained by

(i) Multiplying the above given derived Jacobi conditions by corresponding multiplier from $\{f, g, f^2, fg, g^2 \cdots\}$ and then substituting for $\{f, g, f^2, fg, g^2 \cdots\}$ in terms of $\{x, y\}$ and then
(ii) Adding these expressions, and
(iii) Collecting the coefficients of $\{x, y, x^2, xy, y^2 \cdots\}$, and
(iv) Setting them to zero.

Thus, e.g collecting the coefficients of $x^2$, we get condition like

$$
\begin{aligned}
&(2c_{20}a_{20}+c_{11}b_{20})d_{01}+(d_{11}a_{20}+2d_{02}b_{20})c_{10} \\
&-(2d_{20}a_{20}+d_{11}b_{20})c_{01}-(c_{11}a_{20}+2c_{02}b_{20})d_{10}+\cdots=0
\end{aligned}
\quad --- \text{(2.6.1)}
$$

$$\vdots$$

**3. An Algorithm to Construct Inverse Function Pair:** We now begin with the main part for the two variable case. Our main objective here is to determine the coefficients of the inverse function in terms of the coefficients of the originally given polynomials and to show that only finitely many of these inverse coefficients are actually nonzero. Our main source for completing this task is

(i) The matrix equation (2.2.1) which actually contains four equations to be obtained by equating the respective matrix elements, and

(ii) The principle as well as derived Jacobi conditions.

We see that **these two things are actually enough.**

The four equations contained in the matrix equation (2.2.1), obtained by equating the respective matrix elements, are

$$x_f = g_y \qquad (3.1)$$

$$x_g = -f_y \qquad (3.2)$$

$$y_f = -g_x \qquad (3.3)$$

$$y_g = f_x \qquad (3.4)$$

By expressing each of $\{f, g, x, y\}$ as sum of homogeneous blocks of same total degree we can express these equations as

$$
\begin{aligned}
&\sum_{k=0}^{\infty}\sum_{r=0}^{k}((k+1)-r)c_{((k+1)-r)r}\, f^{(k-r)}g^{r} \\
&=\sum_{k=0}^{l_1}\ \sum_{r=0}^{k}(r+1)b_{(k-r)(r+1)}x^{(k-r)}y^{r}
\end{aligned}
\qquad (3.5)
$$

$$\sum_{k=0}^{\infty} \sum_{r=0}^{k}(r+1)c_{(k-r)(r+1)} f^{(k-r)} g^{r} = -\sum_{k=0}^{l_1} \sum_{r=0}^{k} (r+1)a_{(k-r)(r+1)} x^{(k-r)} y^{r} \tag{3.6}$$

$$\sum_{k=0}^{\infty} \sum_{r=0}^{k} ((k+1)-r)d_{((k+1)-r)r} f^{(k-r)} g^{r} = -\sum_{k=0}^{l_1} \sum_{r=0}^{k} ((k+1)-r)b_{((k+1)-r)r} x^{(k-r)} y^{r} \tag{3.7}$$

$$\sum_{k=0}^{\infty} \sum_{r=0}^{k}(r+1)d_{(k-r)(r+1)} f^{(k-r)} g^{r} = \sum_{k=0}^{l_1} \sum_{r=0}^{k} ((k+1)-r)a_{((k+1)-r)r} x^{(k-r)} y^{r} \tag{3.8}$$

Now, if we substitute the polynomial expressions for $f$ and $g$ in terms of $x$ and $y$ in the above equations and rearrange the resulting equations such that the terms involving the coefficients of only linear parts of $f$ and $g$ and their products, powers etc., i.e. the terms involving **the power products of the coefficients of linear parts of** $f$ and $g$, are **only the things** we keep on the left side of the equations and **all the rest on the right side,** we have the following type of equations relating the coefficients $\{c^s, d^s\}$ to $\{a^s, b^s\}$.

Note that we take the value of the binomial coefficient $\binom{n}{m} = 0$, when $m > n$, and, otherwise, $\binom{n}{m} = \frac{n!}{m!(n-m)!}$ when $m \le n$, as usual.

Each **homogeneous block of degree** $k$, is made up of terms (monomials) $\{x^k, x^{k-1}y, \cdots, xy^{k-1}, y^k\}$.

(i) For a term in this block of the form $x^{k-s}y^s$ we get from equation (3.5)

$$\sum_{r=0}^{k}((k+1)-r)c_{((k+1)-r)r}\sum_{p=0}^{r}\binom{k-r}{s-p}\binom{r}{p}a_{10}^{(k-r-s+p)}a_{01}^{(s-p)}b_{10}^{(r-p)}b_{01}^{p} \quad (3.9)$$

= {Coefficient products arrived at from nonlinear parts of the $f$ and $g$ }

(ii) For a term in this block of the form $x^{k-s}y^{s}$ we get from equation (3.6)

$$\sum_{r=0}^{k}((k+1)-r)d_{((k+1)-r)r}\sum_{p=0}^{r}\binom{k-r}{s-p}\binom{r}{p}a_{10}^{(k-r-s+p)}a_{01}^{(s-p)}b_{10}^{(r-p)}b_{01}^{p} \quad (3.10)$$

= {Coefficient products arrived at from nonlinear parts of the $f$ and $g$ }

(iii) For a term in this block of the form $x^{k-s}y^{s}$ we get from equation (3.7)

$$\sum_{r=0}^{k}(r+1)c_{(k-r)(r+1)}\sum_{p=0}^{r}\binom{k-r}{s-p}\binom{r}{p}a_{10}^{(k-r-s+p)}a_{01}^{(s-p)}b_{10}^{(r-p)}b_{01}^{p} \quad (3.11)$$

= {Coefficient products arrived at from nonlinear parts of the $f$ and $g$ }

(iv) For a term in this block of the form $x^{k-s}y^{s}$ we get from equation (3.8)

$$\sum_{r=0}^{k}(r+1)d_{(k-r)(r+1)}\sum_{p=0}^{r}\binom{k-r}{s-p}\binom{r}{p}a_{10}^{(k-r-s+p)}a_{01}^{(s-p)}b_{10}^{(r-p)}b_{01}^{p} \quad (3.12)$$

= {Coefficient products arrived at from nonlinear parts of the $f$ and $g$ }

In order to understand our method we begin from zero degree block and increase the degree in steps to consider in succession the corresponding higher degree blocks that result with the increase in the total degree in steps.

**3.1 The System of Equations:**

(1) From the block of **total degree = 0**, we have from equations (3.9) to (3.12) (which are actually equations (3.1) to (3.4) with terms partitioned into homogeneous blocks)

$$c_{10} = b_{01} \qquad (3.1.1)$$

$$c_{01} = -a_{01} \qquad (3.1.2)$$

$$d_{10} = -b_{10} \qquad (3.1.3)$$

$$d_{01} = a_{10} \qquad (3.1.4)$$

(1) Similarly, from the block of **total degree = 1** in $x, y$ we have

$$2c_{20}a_{10} + c_{11}b_{10} = b_{11} \qquad (3.1.5)$$

$$2c_{20}a_{01} + c_{11}b_{01} = 2b_{02} \qquad (3.1.6)$$

$$2c_{02}b_{10} + c_{11}a_{10} = -a_{11} \qquad (3.1.7)$$

$$2c_{02}b_{01} + c_{11}a_{01} = -2a_{02} \qquad (3.1.8)$$

$$2d_{20}a_{10} + d_{11}b_{10} = -2b_{20} \qquad (3.1.9)$$

$$2d_{20}a_{01} + d_{11}b_{01} = -b_{11} \qquad (3.1.10)$$

$$2d_{02}b_{10} + d_{11}a_{10} = 2a_{20} \qquad (3.1.11)$$

$$2d_{02}b_{01} + d_{11}a_{01} = a_{11} \qquad (3.1.12)$$

These equations when expressed in matrix form become

$$\begin{pmatrix} a_{10} & b_{10} \\ a_{01} & b_{01} \end{pmatrix} \begin{pmatrix} 2c_{20} \\ c_{11} \end{pmatrix} = \begin{pmatrix} b_{11} \\ 2b_{20} \end{pmatrix} \qquad (3.1.13)$$

$$\begin{pmatrix} a_{10} & b_{10} \\ a_{01} & b_{01} \end{pmatrix} \begin{pmatrix} c_{11} \\ 2c_{02} \end{pmatrix} = \begin{pmatrix} -a_{11} \\ -2a_{20} \end{pmatrix} \qquad (3.1.14)$$

$$\begin{pmatrix} a_{10} & b_{10} \\ a_{01} & b_{01} \end{pmatrix} \begin{pmatrix} 2d_{20} \\ d_{11} \end{pmatrix} = \begin{pmatrix} -2b_{20} \\ -b_{11} \end{pmatrix} \qquad (3.1.15)$$

$$\begin{pmatrix} a_{10} & b_{10} \\ a_{01} & b_{01} \end{pmatrix} \begin{pmatrix} d_{11} \\ 2d_{02} \end{pmatrix} = \begin{pmatrix} 2a_{20} \\ a_{11} \end{pmatrix} \qquad (3.1.16)$$

Let us denote the (2×2)-matrix in the above matrix equations by $M$ . It is clear from the principle Jacobi condition that $\det(M) = 1$. And so, we can solve these invertible systems of matrix equations and determine the values of coefficients $\{c_{20}, c_{11}, c_{02}, d_{20}, d_{11}, d_{02}\}$ uniquely in terms of the coefficients $\{a_{20}, a_{11}, a_{02}, b_{20}, b_{11}, b_{02}, a_{10}, a_{01}, b_{10}, b_{01}\}$.

(2) Similarly, from the block of **total degree = 2** in $x, y$ we have

$$\begin{pmatrix} a_{10}^2 & a_{10}b_{10} & b_{10}^2 \\ 2a_{10}a_{01} & (a_{10}b_{01} + b_{10}a_{01}) & 2b_{10}b_{01} \\ a_{01}^2 & a_{01}b_{01} & b_{01}^2 \end{pmatrix} \begin{pmatrix} 3c_{30} \\ 2c_{21} \\ c_{12} \end{pmatrix}$$

$$= \begin{pmatrix} b_{21} - (2c_{20}a_{20} + c_{11}b_{20}) \\ 2b_{12} - (2c_{20}a_{11} + c_{11}b_{11}) \\ 3b_{03} - (2c_{20}a_{02} + c_{11}b_{02}) \end{pmatrix} \qquad \cdots (3.1.17)$$

$$\begin{pmatrix} a_{10}^2 & a_{10}b_{10} & b_{10}^2 \\ 2a_{10}a_{01} & (a_{10}b_{01} + b_{10}a_{01}) & 2b_{10}b_{01} \\ a_{01}^2 & a_{01}b_{01} & b_{01}^2 \end{pmatrix} \begin{pmatrix} c_{21} \\ 2c_{12} \\ 3c_{03} \end{pmatrix}$$

$$= -\begin{pmatrix} a_{21} + (2c_{02}b_{20} + c_{11}a_{20}) \\ 2a_{12} + (2c_{02}b_{11} + c_{11}a_{11}) \\ 3a_{03} + (2c_{02}b_{02} + c_{11}a_{02}) \end{pmatrix} \qquad \cdots (3.1.18)$$

$$\begin{pmatrix} a_{10}^2 & a_{10}b_{10} & b_{10}^2 \\ 2a_{10}a_{01} & (a_{10}b_{01} + b_{10}a_{01}) & 2b_{10}b_{01} \\ a_{01}^2 & a_{01}b_{01} & b_{01}^2 \end{pmatrix} \begin{pmatrix} 3d_{30} \\ 2d_{21} \\ d_{12} \end{pmatrix}$$

$$= -\begin{pmatrix} 3b_{30} + (2d_{20}a_{20} + d_{11}b_{20}) \\ 2b_{21} + (2d_{20}a_{11} + d_{11}b_{11}) \\ b_{12} + (2d_{20}a_{02} + d_{11}b_{02}) \end{pmatrix} \qquad \cdots (3.119)$$

$$\begin{pmatrix} a_{10}^2 & a_{10}b_{10} & b_{10}^2 \\ 2a_{10}a_{01} & (a_{10}b_{01} + b_{10}a_{01}) & 2b_{10}b_{01} \\ a_{01}^2 & a_{01}b_{01} & b_{01}^2 \end{pmatrix} \begin{pmatrix} d_{21} \\ 2d_{12} \\ 3d_{03} \end{pmatrix}$$

$$= \begin{pmatrix} 3a_{30} - (2d_{02}b_{20} + d_{11}a_{20}) \\ 2a_{21} - (2d_{02}b_{11} + d_{11}a_{11}) \\ a_{12} - (2d_{02}b_{02} + d_{11}a_{02}) \end{pmatrix} \quad \cdots\cdot(3.120)$$

**Remark 3.1.1:** Note that the (3×3)-matrices in the above four equations are identical and each one is equal to transpose of the matrix $U_3$(defined in the proof of the theorem 2.3.1), i.e. equal to $U_3^T$ . Since matrix $U_3$ is invertible and $(U_3)^{-1} = V_3$, therefore, $U_3^T$ will be invertible and $(U_3^T)^{-1} = V_3^T$ .
$(V_3 = (U_3)^{-1} \Leftrightarrow U_3V_3 = I_3 \Leftrightarrow V_3^T U_3^T = I_3^T = I_3 \Leftrightarrow V_3^T = (U_3^T)^{-1})$.
Thus, we can solve these invertible systems of matrix equations and determine the values of coefficients $\{c_{30}, c_{21}, c_{12}, c_{03}, d_{30}, d_{21}, d_{12}, d_{02}\}$ uniquely in terms of the coefficients
$\{a_{30}, a_{21}, a_{12}, a_{03}, b_{30}, b_{21}, b_{12}, b_{03}, a_{20}, a_{11}, a_{02}, b_{20}, b_{11}, b_{02},$
$a_{10}, a_{01}, b_{10}, b_{01}\}$.

**Remark 3.1.2:** If we proceed on similar lines for the block of **total degree = (n−1)**, we will have from equations (3.9) to (3.12) (which are actually equations (3.1) to (3.4) with terms partitioned into homogeneous blocks) four invertible matrix equations, namely,

$$U_n^T \begin{pmatrix} nc_{n0} \\ (n-1)c_{(n-1)1} \\ (n-2)c_{(n-2)2} \\ \vdots \\ 2c_{2(n-2)} \\ c_{1(n-1)} \end{pmatrix} = \begin{pmatrix} b_{(n-1)1} - Y_1^1 \\ 2b_{(n-2)2} - Y_2^1 \\ 3b_{(n-3)3} - Y_3^1 \\ \vdots \\ (n-1)b_{1(n-1)} - Y_{(n-1)}^1 \\ nb_{0n} - Y_n^1 \end{pmatrix} \quad \cdots\cdot(3.1.21)$$

$$U_n^T \begin{pmatrix} c_{(n-1)1} \\ 2c_{(n-2)2} \\ 3c_{(n-3)3} \\ \vdots \\ (n-1)c_{1(n-1)} \\ nc_{0n} \end{pmatrix} = -\begin{pmatrix} a_{(n-1)1} + Y_1^2 \\ 2a_{(n-2)2} + Y_2^2 \\ 3a_{(n-3)3} + Y_3^2 \\ \vdots \\ (n-1)a_{1(n-1)} + Y_{(n-1)}^2 \\ na_{0n} + Y_n^2 \end{pmatrix} \quad \cdots\cdots (3.1.22)$$

$$U_n^T \begin{pmatrix} nd_{n0} \\ (n-1)d_{(n-1)1} \\ (n-1)d_{(n-2)2} \\ \vdots \\ 2d_{2(n-2)} \\ d_{1(n-1)} \end{pmatrix} = -\begin{pmatrix} nb_{n0} + Y_1^3 \\ (n-1)b_{(n-1)1} + Y_2^3 \\ (n-2)b_{(n-2)2} + Y_3^3 \\ \vdots \\ 2b_{2(n-2)} + Y_{(n-1)}^3 \\ b_{1(n-1)} + Y_n^3 \end{pmatrix} \quad \cdots\cdots (3.1.23)$$

$$U_n^T \begin{pmatrix} d_{(n-1)1} \\ 2d_{(n-2)2} \\ 3d_{(n-3)3} \\ \vdots \\ (n-1)d_{1(n-1)} \\ nd_{0n} \end{pmatrix} = \begin{pmatrix} na_{n0} - Y_1^4 \\ (n-1)a_{(n-1)1} - Y_2^4 \\ (n-2)a_{(n-2)2} - Y_3^4 \\ \vdots \\ 2a_{2(n-2)} - Y_{(n-1)}^4 \\ a_{1(n-1)} - Y_n^4 \end{pmatrix} \quad \cdots\cdots (3.1.24)$$

Thus, we can solve these invertible systems of matrix equations and determine the values of coefficients $\{c_{n0},\cdots,c_{0n},d_{n0},\cdots,d_{0n}\}$ uniquely in terms of the coefficients

$\{a_{n0},\cdots,a_{0n},b_{n0},\cdots,b_{0n},a_{(n-1)0},\cdots,a_{0(n-1)},b_{(n-1)0},\cdots$

$,b_{0(n-1)},\cdots,a_{10},a_{01},b_{10},b_{01}\}$.

Let us call the terms collected in each of $Y_i^j, i=1,2,\cdots,$ and $j=1,2,3,4$, the **"residuum"**.

**Remark 3.1.2:** Since $f$ and $g$ are polynomials their coefficients that appear on the left hand side of the residuum in the vectors on the right hand side of the above equations (3.2.21) to (3.2.24) will subsequently evaporate (vanish as they being zero after a certain stage) and thus, settling the plane Jacobian conjecture thus reduces to showing that the **residuum also vanishes** after that stage of the evaporation of the corresponding $a^s$ and $b^s$. We will see that the vanishing of the residuum follows from the derived Jacobi conditions!!

Our procedure to deal with the Jacobian problem is **exactly same** for each case. We now proceed to illustrate this procedure we begin (because of its simplicity) with the example of polynomials of degree two and (in two variables):

**Example 3.1.1:** Let $f$ and $g$ be polynomials with quadratic total degree. Thus,

$$f = a_{10}x + a_{01}y + a_{20}x^2 + a_{11}xy + a_{02}y^2$$
$$g = b_{10}x + b_{01}y + b_{20}x^2 + b_{11}xy + b_{02}y^2$$

with coefficients in the field $k$ of characteristic zero and the Jacobian

$$J_{(x,y)}(f,g) = \det\begin{pmatrix} f_x & f_y \\ g_x & g_y \end{pmatrix} = 1 \neq 0$$

Proceeding as is done in the section 3.1 we can obtain equations same as equations (3.1.1) to (3.1.20). Solving matrix equations (3.1.13) to (3.1.16) we can obtain $c_{20}, c_{11}, c_{02}$ as follows:

$$2c_{20} = b_{01}b_{11} - 2b_{10}b_{02}$$
$$c_{11} = -a_{01}b_{11} + 2a_{10}b_{02} \text{, also,}$$
$$c_{11} = -b_{01}a_{11} + 2b_{10}a_{02}$$
$$2c_{02} = a_{01}a_{11} - 2a_{10}a_{02}$$

Since $f$ and $g$ are quadratic polynomials, in equations (3.1.17) to (3.1.20) actually the coefficients $\{a_{30}, a_{21}, a_{12}, a_{03}, b_{30}, b_{21}, b_{12}, b_{03}\}$ are all equal to zero. So, we put zero value for them in these equations. We then show that the residuum in equation (3.1.17) vanishes. Now, after substitution of zero value as mentioned just above for cubic coefficients and after substitution for $c_{20}, c_{11}, c_{02}$ from the above equations the vector on the right side of equation (3.1.17) becomes

$$\begin{pmatrix} 2c_{20}a_{20} + c_{11}b_{20} \\ 2c_{20}a_{11} + c_{11}b_{11} \\ 2c_{20}a_{02} + c_{11}b_{02} \end{pmatrix} = \begin{pmatrix} b_{01}(a_{20}b_{11} - b_{20}a_{11}) + 2b_{10}(b_{20}a_{02} - a_{20}b_{02}) \\ b_{01}(a_{11}b_{11} - b_{11}a_{11}) + 2b_{10}(b_{11}a_{02} - a_{11}b_{02}) \\ b_{01}(b_{11}a_{02} - a_{11}b_{02}) + 2b_{10}(b_{02}a_{02} - a_{02}b_{02}) \end{pmatrix}$$

Using the derived Jacobian conditions for $\{x^2, xy, y^2\}$ given in the equations (2.4.4) to (2.4.6) it is clear that this vector (let us call it the **residuum vector**) is actually equal to zero!!

By proceeding on similar lines one can show that the corresponding residuum vector in equations (3.1.18) to (3.1.20) also vanishes!!!

Thus, we get the inverse functions which are actually the following polynomials:

$$x = b_{01} f - a_{01} g + \frac{1}{2}(b_{01}b_{11} - 2b_{10}b_{02}) f^2 + (2a_{10}b_{02} - a_{01}b_{11}) fg$$
$$+ \frac{1}{2}(a_{01}a_{11} - 2a_{10}a_{02}) g^2$$

$$y = -b_{10} f + a_{10} g + \frac{1}{2}(b_{10}b_{11} - 2b_{01}b_{20}) f^2 + (2a_{01}b_{20} - a_{10}b_{11}) fg$$
$$+ \frac{1}{2}(a_{10}a_{11} - 2a_{01}a_{20}) g^2$$

**Remark 3.1.1:** Note that the terms, like e.g. $2c_{20}a_{20} + c_{11}b_{20}$, that appear in the residuum vector also appear as terms in the mixed Jacobi conditions (refer to the mixed Jacobi condition for $x^2$, namely, equation (2.6.1)) and they vanish as just seen.

**Example 3.1.2:** Let $f$ and $g$ be polynomials with cubic total degree. Thus,

$$f = a_{10}x + a_{01}y + a_{20}x^2 + a_{11}xy + a_{02}y^2 + a_{30}x^3 + a_{21}x^2y + a_{12}xy^2 + a_{03}y^3$$

$$g = b_{10}x + b_{01}y + b_{20}x^2 + b_{11}xy + b_{02}y^2 + b_{30}x^3 + b_{21}x^2y + b_{12}xy^2 + b_{03}y^3$$

with coefficients in the field $k$ of characteristic zero and the Jacobian

$$J_{(x,y)}(f,g) = \det\begin{pmatrix} f_x & f_y \\ g_x & g_y \end{pmatrix} = 1 \neq 0$$

We proceed as is done in the above example and continue further with determining the values of coefficients $\{c_{30}, c_{21}, c_{12}, c_{03}, d_{30}, d_{21}, d_{12}, d_{02}\}$ uniquely in terms of the coefficients

$$\{a_{30}, a_{21}, a_{12}, a_{03}, b_{30}, b_{21}, b_{12}, b_{03}, a_{20}, a_{11}, a_{02}, b_{20}, b_{11}, b_{02}, \\ a_{10}, a_{01}, b_{10}, b_{01}\}.$$

We then form the next four equations, the first one of which is

$$U_4^T \begin{pmatrix} 4c_{40} \\ 3c_{31} \\ 2c_{22} \\ c_{13} \end{pmatrix} = \begin{pmatrix} b_{31} - Y_1 \\ 2b_{22} - Y_2 \\ 3b_{13} - Y_3 \\ 4b_{04} - Y_4 \end{pmatrix}$$

where

$$Y_1 = \{3c_{30}\, 2a_{10}a_{20} + 2c_{21}(a_{20}b_{10} + b_{20}a_{10}) + c_{12}\, 2b_{10}b_{20}\}$$

$$Y_2 = \{3c_{30}(2a_{20}a_{01} + 2a_{11}a_{10}) + 2c_{21}(a_{20}b_{01} + b_{20}a_{01} + \\ a_{11}b_{10} + b_{11}a_{10}) + c_{12}(2b_{20}b_{01} + 2b_{11}b_{10})\}$$

$$Y_3 = \{3c_{30}(2a_{02}a_{10} + 2a_{11}a_{01}) + 2c_{21}(a_{02}b_{10} + b_{02}a_{10} + \\ a_{11}b_{01} + b_{11}a_{01}) + c_{12}(2b_{02}b_{01} + 2b_{11}b_{01})\}$$

$$Y_4 = \{3c_{30}\, 2a_{01}a_{02} + 2c_{21}(a_{02}b_{01} + b_{02}a_{01}) + c_{12}\, 2b_{01}b_{02}\}$$

Since $f, g$ are cubic, coefficients of the monomials in the homogeneous block of total degree = 4 will have zero value. Therefore, in the right hand side vector in the above equation, $b_{31} = b_{22} = b_{13} = b_{04} = 0$. Also, one can check using the derived Jacobi conditions that actually $Y_i, i = 1,2,3,4$, **are all equal to zero!!** Thus, all the components of the vector on the left side of the above equation vanish! Proceeding with other three equations one can show that actually all the fourth (and higher) degree coefficients $c_{ij}, d_{ij}$ where $i + j \geq 4$ vanish! So, finally one obtains the inverse functions (as is done above for quadratic case) as actually the cubic polynomials!! This example implies that we have verified Jacobian conjecture for a cubic polynomial of most general form (and not only for the special form which is actually enough to be checked by the important degree reduction theorem due to H. Bass, E. H. Connell, and D. Wright [3]) for two variable case.

**Remark 3.1.3: Fortunately, we are empowered** by result [3] about the reduction in degree. This result very much simplifies our job of checking the vanishing of the terms $Y_i, i = 1,2,3,4$ forming the so called **residuum**. So we now proceed to demonstrate its use here:
As per the result of BCW in [3], we can take

$$f = a_{10}x + a_{01}y + a_{30}x^3 + a_{21}x^2y + a_{12}xy^2 + a_{03}y^3$$

$$g = b_{10}x + b_{01}y + b_{30}x^3 + b_{21}x^2y + b_{12}xy^2 + b_{03}y^3$$

and proceed in the same steps as is done in example 3.1.2 above and obtain $Y_i, i = 1,2,3,4$. Now it is clear to see that each term in each $Y_i, i = 1,2,3,4$ contains **a coefficient of quadratic term,** i.e. some one of the coefficients among $\{a_{20}, a_{11}, a_{02}, b_{20}, b_{11}, b_{02}\}$**,** which is already **zero due to the form of the above polynomials** $f, g$ **in accordance with BCW result!!** Thus, for such polynomials all $Y_i, i = 1,2,3,4$ are **automatically zero!!!** Thus, by using the **BCW form** for the polynomials the **Residuum vanishes automatically** as desired**.** Also, it is **important** to note that **when given polynomial have BCW form** then quadratic terms are absent in the given polynomials, therefore, all coefficients of quadratic terms are equal to zero, i.e. all $a_{ij} = 0, b_{ij} = 0$ whenever $i + j = 2$ and this further implies from equations (3.1.13) to (3.1.16) that also $c_{ij} = 0, d_{ij} = 0$ whenever $i + j = 2$ **when given polynomials have BCW form.**

**Remark 3.1.4:** For the case of quadratic polynomials, $f$ and $g$, having nonzero Jacobian we have obtained the inverse functions $x$ and $y$ as given above. Now the natural question to be asked and answered is what will be the inverse functions $x$ and $y$ when $f$ and $g$ will be cubic polynomials of BCW form having nonzero Jacobian (=1)? We now proceed to answer this question. Thus, we can take

$$f = a_{10}x + a_{01}y + a_{30}x^3 + a_{21}x^2 y + a_{12}xy^2 + a_{03}y^3$$

$$g = b_{10}x + b_{01}y + b_{30}x^3 + b_{21}x^2 y + b_{12}xy^2 + b_{03}y^3$$

Using the development done in Example 3.1.2 we can grant that for the present case the inverse functions will be actually the **cubic polynomials** as given below:

$$x = \sum_{m=0}^{3} \sum_{n=0}^{3} c_{mn} f^m g^n \qquad (3.1.25)$$

and

$$y = \sum_{m=0}^{3} \sum_{n=0}^{3} d_{mn} f^m g^n \qquad (3.1.26)$$

$c_{mn}, d_{mn} \in k$, and $m + n \le 3$

Now, using the theory developed so far we can now easily construct the inverse functions $x$ and $y$ when $f$ and $g$ will be cubic polynomials of BCW form as given above having nonzero Jacobian (= 1). These inverse functions are actually cubic polynomials as given below:

$$\begin{aligned} x &= b_{01}f - a_{01}g + \frac{1}{3}(b_{01}^2 b_{21} - 2b_{01}b_{10}b_{12} + 3b_{10}^2 b_{03}) f^3 \\ &+ (-b_{01}^2 a_{21} + 2b_{01}b_{10}a_{12} - 3b_{10}^2 a_{03}) f^2 g \\ &+ (a_{01}^2 b_{21} - 2a_{01}a_{10}b_{12} + 3a_{10}^2 b_{03}) fg^2 \\ &+ \frac{1}{3}(-a_{01}^2 a_{21} + 2a_{01}a_{10}a_{12} - 3a_{10}^2 a_{03}) g^3 \end{aligned} \qquad (3.1.27)$$

$$
\begin{aligned}
y &= -b_{10} f + a_{10} g + \frac{1}{3}(-3b_{01}^2 b_{30} + 2b_{01}b_{10}b_{21} - b_{10}^2 b_{12}) f^3 \\
&+ (3b_{01}^2 a_{30} - 2b_{01}b_{10}a_{21} + b_{10}^2 a_{12}) f^2 g \\
&+ (-3a_{01}^2 b_{30} + 2a_{01}a_{10}b_{21} - a_{10}^2 b_{12}) fg^2 \\
&+ \frac{1}{3}(3a_{01}^2 a_{30} - 2a_{01}a_{10}a_{21} + a_{10}^2 a_{12}) g^3
\end{aligned} \tag{3.1.28}
$$

We can proceed along almost same lines for considering the case when $f$ and $g$ will be cubic polynomials of BCW form as given above having nonzero Jacobian (= 1) and the only change is the change of dimension, i.e. when we will be dealing with the case with number of variables $n = 3,4,\cdots,k,\cdots$ instead of the above considered case of $n = 2$. By following exactly the same steps with appropriate changes required by the change in the dimension we can easily obtain inverse functions which will be cubic polynomials of BCW form as desired. Thus, the solution for the Jacobian question!!

**4. More about the Matrices $U_n^T$ and $V_n^T$ :** The nonsingular nature of the matrices $U_n^T$ for all $n$ can also be directly verified. This will be done in this section. We show that the value of the determinant of these matrices is equal to unity. In fact, these determinants are powers of Jacobian, and these powers form a line in the Pascal triangle.

Using equations (3.1.17) to (3.1.20) we can proceed with the construction of the matrices $U_n^T$ corresponding to higher and higher degree homogeneous blocks. But they can also be obtained by defining a simple product rule for lower size matrices. We will describe such a rule in this section.

We state the formulae for inverse functions in the case when $f$ is a polynomial of total degree $k$ and $g$ is linear which can be straightforwardly obtained.

**4.1. The nonsingular nature of the Matrices $U_n^T$ and $V_n^T$ :**

We put $\left(\dfrac{a_{10}}{b_{10}}\right) = u$ and $\left(\dfrac{a_{01}}{b_{01}}\right) = v$. Therefore, the Jacobian will be

$$(b_{10}b_{01})\det\begin{pmatrix} u & 1 \\ v & 1 \end{pmatrix} = (b_{10}b_{01})(u-v) = M = 1 \neq 0 .$$

Now, let us denote the Jacobian by $M$ and the (3×3) matrix on the left hand side of equations (3.1.17) to (3.1.20) by $U_3^T$, then we have

$$\det(U_3^T) = (b_{10}b_{01})^3 \det\begin{pmatrix} u^2 & u & 1 \\ 2uv & (u+v) & 2 \\ v^2 & v & 1 \end{pmatrix}$$

Now, on the matrix inside the determinant operator we perform the following actions:

(i) We multiply the last column by $v$ and subtract it from second column.

(ii) We multiply the second column by $v$ and subtract it from first column

This leads to

$$\det(U_3^T) = (b_{10}b_{01})^3 (u-v)^2 \det\begin{pmatrix} u & 1 & 1 \\ v & 1 & 2 \\ 1 & 0 & 1 \end{pmatrix} = (b_{10}b_{01})^3 (u-v)^3$$

Thus,

$$\det(U_3^T) = M^3 = 1 \neq 0 .$$

If we proceed with equations like equations (3.1.17) to (3.1.20) for the next set and form, say, $U_4^T$ and carry out similar procedures as above we will see that we get

$$\det(U_4^T) = (b_{10}b_{01})^6 \det\begin{pmatrix} u^3 & u^2 & u & 1 \\ 3u^2v & u^2+2uv & 2u+v & 3 \\ 3uv^2 & v^2+2uv & 2v+u & 3 \\ v^3 & v^2 & v & 1 \end{pmatrix}$$

therefore,

$$\det(U_4^T) = (b_{10}b_{01})^6 (u-v)^3 \det\begin{pmatrix} u^2 & u & 1 & 1 \\ 2uv & (u+v) & 2 & 3 \\ v^2 & v & 1 & 3 \\ 0 & 0 & 0 & 1 \end{pmatrix}$$

So, from the previous stage,

$$\det(U_4^T) = (b_{10}b_{01})^6 (u-v)^6 = M^6 = 1 \neq 0$$

We can now prove the following

**Theorem 4.1:** The determinant of the matrix $U_n^T$, $n = 2,3,\cdots$ is power of Jacobian $M$, i.e. $\det(U_n^T) = M^k$ for some positive integer $k$. Moreover, the values $k$ considered in succession form a line (**third diagonal line**) in the Pascal triangle. Thus,

| $n$ (size) | 2 | 3 | 4 | 5 | 6 | $\cdots$ |
|---|---|---|---|---|---|---|
| $k$ (power) | 1 | 3 | 6 | 10 | 15 | $\cdots$ |

**Proof:** It follows from simple induction.

**Step 1:** The result is true for $n = 2,3,4$ as seen above.

**Step 2:** We assume the result for $n = (k-1)$ and prove the it for $n = k$.

Using any of the equations (3.9) to (3.12) we can construct matrix $U_k^T$. By substituting $\left(\frac{a_{10}}{b_{10}}\right) = u$ and $\left(\frac{a_{01}}{b_{01}}\right) = v$, as above, we can transform the elements of this matrix in terms of $u$ and $v$. Now, if we multiply $(k-i)$-th column of this matrix by $v$ and subtract it from $(k-i-1)$-th column for each $i = 0,\cdots,(k-2)$ we get the determinant as

$$\det(U_k^T) = const \times (u-v)^k \det\begin{pmatrix} & & & \cdot \\ & B & & \vdots \\ & & & \cdot \\ 0 & \cdots & 0 & 1 \end{pmatrix}$$

where $B$ is the matrix in terms of $u, v$ corresponding to previous i.e. $(k-1)$-th stage. Hence etc. □

**4.2 A Special Case for Inverse Functions for $f$ with total degree $k$ and a linear $g$:** Let $f$ be of total degree $k$ and let $g$ be linear i.e.

$$f = a_{10}x + a_{01}y + a_{20}x^2 + a_{11}xy + a_{02}y^2 + \cdots + a_{0k}y^k$$

and

$$g = b_{10}x + b_{01}y$$

such that

$$J_{(x,y)}(f,g) = \det\begin{pmatrix} f_x & f_y \\ g_x & g_y \end{pmatrix} = 1 \neq 0$$

In this case, it is an easy exercise to check that we (can easily) get the inverse functions as

$$x = b_{01} f - a_{01} g - \frac{a_{02}}{b_{01}} g^2 - \frac{a_{03}}{b_{01}^2} g^3 - \cdots - \frac{a_{0k}}{b_{01}^{(k-1)}} g^k$$

$$y = -b_{10} f + a_{10} g + \frac{a_{11}}{2b_{01}} g^2 + \frac{a_{12}}{3b_{01}^2} g^3 + \cdots + \frac{a_{1(k-1)}}{kb_{01}^{(k-1)}} g^k$$

**5. The Direct or Algorithmic Approach for Several Variables:** In this section we will briefly see how the considerations developed so far for the two variables can be extended for the several variables case by proceeding on similar lines. To illustrate the approach we first consider the **three variable case** in more detail so that it will make absolutely clear how the same approach can take us smoothly and without any hindrance and doubt into several variable case.

**5.1 The Three Variables Case:** Suppose we are given the following polynomials for which the Jacobian is a nonzero constant (=1) as desired.

f=a100*x+a010*y+a001*z+a200*x^2+a110*x*y+a101*x*z+a020*y^2+a011*y*z+a002*z^2+a300*x^3+a210*x^2*y+a201*x^2*z+a120*x*y^2+a111*x*y*z+a102*x*z^2+a030*y^3+a021*y^2*z+a012*y*z^2+a003*z^3

(5.1.1)

and

g=b100*x+b010*y+b001*z+b200*x^2+b110*x*y+b101*x*z+b020*y^2+b011*y*z+b002*z^2+b300*x^3+b210*x^2*y+b201*x^2*z+b120*x*y^2+b111*x*y*z+b102*x*z^2+b030*y^3+b021*y^2*z+b012*y*z^2+b003*z^3 (5.1.2)

and

h=c100*x+c010*y+c001*z+c200*x^2+c110*x*y+c101*x*z+c020*y^2+c011*y*z+c002*z^2+c300*x^3+c210*x^2*y+c201*x^2*z+c120*x*y^2+c111*x*y*z+c102*x*z^2+c030*y^3+c021*y^2*z+c012*y*z^2+c003*z^3

aijk, bijk, cijk $\in k$, and actually aijk = 0, bijk = 0, cijk = 0 when i+j+k = 2, i.e. f, g, h have BCW form. Also, the Jacobian

$$J_{(x,y,z)}(f,g,h) = \det\begin{bmatrix} f_x & f_y & f_z \\ g_x & g_y & g_z \\ h_x & h_y & h_z \end{bmatrix} = 1 \neq 0 \qquad (5.1.3)$$

We need to show that $x, y, z$ (which we can express as power series in $f, g, h$ using the inverse function theorem) are **actually polynomials**, i.e.

x=p100*f+p010*g+p001*h+p200*f^2+p110*f*g+p101*f*h+p020*g^2+p011*g*h+p002*h^2+p300*f^3+p210*f^2*g+p201*f^2*h+p120*f*g^2+p111*f*g*h+p102*f*h^2+p030*g^3+p021*g^2*h+p012*g*h^2+p003*h^3 (5.1.4)

and

y=q100*f+q010*g+q001*h+q200*f^2+q110*f*g+q101*f*h+q020*g^2+q011*g*h+q002*h^2+q300*f^3+q210*f^2*g+q201*f^2*h+q120*f*g^2+q111*f*g*h+q102*f*h^2+q030*g^3+q021*g^2*h+q012*g*h^2+q003*h^3 (5.1.5)

and

z=r100*f+r010*g+r001*h+r200*f^2+r110*f*g+r101*f*h+r020*g^2+r011*g*h+r002*h^2+r300*f^3+r210*f^2*g+r201*f^2*h+r120*f*g^2+r111*f*g*h+r102*f*h^2+r030*g^3+r021*g^2*h+r012*g*h^2+r003*h^3 (5.1.6)

where pijk, qijk, rijk $\in k$.

**5.2 The Local Homeomorphism:** Let $F(x,y,z) = (f,g,h)$. Let $f_x, f_y, f_z, g_x, g_y, g_z, h_x, h_y, h_z$ exist and are continuous, therefore,

$$DF(x,y,z)=\begin{bmatrix} f_x & f_y & f_z \\ g_x & g_y & g_z \\ h_x & h_y & h_z \end{bmatrix}$$ exists and since

$$\det(DF(x,y,z)=J_{(x,y,z)}(f,g,h)=1\neq 0\,,$$

therefore, by the inverse function theorem, $F^{-1}$ exists and moreover

$$DF^{-1}(f,g,h)=\left(DF(x,y,z)\right)^{-1}.$$

i.e.

$$\begin{bmatrix} x_f & x_g & x_h \\ y_f & y_g & y_h \\ z_f & z_g & z_h \end{bmatrix}=\begin{bmatrix} f_x & f_y & f_z \\ g_x & g_y & g_z \\ h_x & h_y & h_z \end{bmatrix}^{-1} \qquad (5.2.1)$$

The (local) inverse is made up of power series in general but the Jacobian conjecture demands that they are in fact polynomials (whose coefficients can be uniquely determined in terms of the coefficients of the originally given polynomials $f,g,h$). By evaluating the inverse on the right hand side of the above equation (5.2.1) and using the fact that the value of the Jacobian is constant and equal to unity we get

$$\begin{bmatrix} x_f & x_g & x_h \\ y_f & y_g & y_h \\ z_f & z_g & z_h \end{bmatrix}=\begin{bmatrix} g_yh_z-g_zh_y & f_zh_y-f_yh_z & f_yg_z-f_zg_y \\ g_xh_z-g_zh_x & f_xh_z-f_zh_x & f_zg_x-f_xg_z \\ g_xh_y-g_yh_x & f_yh_x-f_xh_y & f_xg_y-f_yg_x \end{bmatrix} \qquad (5.2.2)$$

**5.3 Monomials in Variables as a Vector Space Basis:** Suppose we are given the following invertible system of equations:

f1=a100*x+a010*y+a001*z (5.3.1)

g1=b100*x+b010*y+b001*z (5.3.2)

h1=c100*x+c010*y+c001*z (5.3.3)

with $\det\begin{pmatrix} a_{100} & a_{010} & a_{001} \\ b_{100} & b_{010} & b_{001} \\ c_{100} & c_{010} & c_{001} \end{pmatrix} = 1 \neq 0$.

Thus, we get

$$\begin{bmatrix} (f1) \\ (g1) \\ (h1) \end{bmatrix} = U \begin{bmatrix} x \\ y \\ z \end{bmatrix}$$

and where

$$U = \begin{bmatrix} a100 & a010 & a001 \\ b100 & b010 & b001 \\ c100 & c010 & c001 \end{bmatrix}$$

It is easy to check that we get the following inverse relations

x=(b010*c001-b001*c010)*f1+(a001*c010-a010*c001)*g1+(a010*b001-a001*b010)*h1 (5.3.4)

y=(-b100*c001+b001*c100)*f1+(a100*c001-a001*c100)*g1+(a100*b001-a001*b100)*h1 (5.3.5)

z=(b100*c010-b010*c100)*f1+(a010*c100-a100*c010)*g1+(a100*b010-a010*b100)*h1 (5.3.6)
Thus, we have

$$\begin{bmatrix} x \\ y \\ z \end{bmatrix} = V \begin{bmatrix} (f1) \\ (g1) \\ (h1) \end{bmatrix}$$

Now, from the process of formation of the above matrix equations (or by actually checking) it is easy to see that $V = U^{-1}$

$$V = \begin{bmatrix} (b010*c001-b001*c010) & (a001*c010-a010*c001) & (a010*b001-a001*b010) \\ (b001*c100-b100*c001) & (a100*c001-a001*c100) & (a100*b001-a001*b100) \\ (b100*c010-b010*c100) & (a010*c100-a100*c010) & (a100*b010-a010*b100) \end{bmatrix}$$

We now state the following theorem which can be proved by proceeding on the similar lines as is done for theorem 2.3.1.

**Theorem 5.3.1:** If the system of equations represented by the equations (5.3.1) to (5.3.3) given above is an invertible system then the system of equations formed by the homogeneous blocks of equations of degree $n$ is also invertible for any $n$.

□

To elaborate the case of degree two we make use of relations (5.3.1) to (5.3.3) and form the following equations (**direct relations)**:

(f1)^2=a100^2*x^2+2*a100*x*a010*y+2*a100*x*a001*z +a010^2*y^2+2*a010*y*a001*z+a001^2*z^2 (5.3.7)

(f1*g1)=a100*x^2*b100+a100*x*b010*y+a100*x*b001*z+a010*y*b100*x+a010*y^2*b010+a010*y*b001*z+a001*z*b100*x +a001*z*b010*y+a001*z^2*b001 (5.3.8)

(f1*h1)=a100*x^2*c100+a100*x*c010*y+a100*x*c001*z+a010*y*c100*x+a010*y^2*c010+a010*y*c001*z+a001*z*c100*x+a001*z*c010*y+a001*z^2*c001 (5.3.9)

(g1)^2=b100^2*x^2+2*b100*x*b010*y+2*b100*x*b001*z+b010^2*y^2+2*b010*y*b001*z+b001^2*z^2 (5.3.10)

(g1*h1)=b100*x^2*c100+b100*x*c010*y+b100*x*c001*z+b010*y*c100*x+b010*y^2*c010+b010*y*c001*z+b001*z*c100*x+b001*z*c010*y+b001*z^2*c001 (5.3.11)

(h1)^2=c100^2*x^2+2*c100*x*c010*y+2*c100*x*c001*z+c010^2*y^2+2*c010*y*c001*z+c001^2*z^2 (5.3.12)

from these equations we can form the following single matrix equation, namely,

$$\begin{bmatrix} (f1)^2 \\ (f1*g1) \\ (f1*h1) \\ (g1)^2 \\ (g1*h1) \\ (h1)^2 \end{bmatrix} = A \begin{bmatrix} (x)^2 \\ (x*y) \\ (x*z) \\ (y)^2 \\ (y*z) \\ (z)^2 \end{bmatrix} \quad (5.3.13)$$

where $A = (sij)$ is a six by six invertible matrix such that $\det(A) = 1$ and the elements of $A$ are as follows:

s11=(a100)^2, s12=2*a100*a010, s13=2*a100*a001, s14=(a010)^2, s15=2*a010*a001, s16=(a001)^2.

s21=a100*b100, s22=(a100*b010+a010*b100), s23=(a100*b001+a001*b100), s24=a010*b010, s25(a010*b001+a001*b010), s26=a001*b001.

s31=a100*c100, s32=(a100*c010+a010*c100),
s33=(a100*c001+a001*c100), s34=a010*c010,
s35=(a010*c001+a001*c010), s36=a001*c001.

s41=(b100)^2, s42=2*b100*b010, s43=2*b100*b001,
s44=(b010)^2, s45=2*b010*b001, s46=(b001)^2.

s51=b100*c100, s52=((b100*c010+b010*c100),
s53=(b100*c001+b001*c100), s54=b010*c010,
s55=(b010*c001+b001*c010), s56=b001*c001.

s61=(c100)^2, s62=2*c100*c010, s63=2*c100c001, s64=(c010)^2,
s65=2*c010*c001, s66=(c001)^2.

Similarly, we make use of relations (5.3.4) to (5.3.6) and form the following equations (**inverse relations)**:

(x)^2=(b010*c001-b001*c010)^2*(f1)^2
+2*(b010*c001-b001*c010)* (-a010*c001+a001*c010)*f1*g1
+2*(b010*c001-b001*c010)* (a010*b001-a001*b010))*f1*h1
+(-a010*c001+a001*c010)^2 *(g1)^2
+2*(-a010*c001+a001*c010)*(a010*b001-a001*b010)*g1*h1
+(a010*b001-a001*b010)^2*(h1)^2 (5.3.14)

x*y=(b010*c001-b001*c010)*(-b100*c001+b001*c100)*(f1)^2
+((b010*c001-b001*c010)*(a100*c001-a001*c100)
+(-a010*c001+a001*c010)*( -b100*c001+b001*c100))*f1*g1
+((b010*c001-b001*c010)* (a100*b001-a001*b100)
+(a010*b001-a001*b010)*( -b100*c001+b001*c100))*f1*h1
+(-a010*c001+a001*c010)*( a100*c001-a001*c100)*(g1)^2
+((-a010*c001+a001*c010)*( a100*b001-a001*b100)
+(a010*b001-a001*b010)*(a100*c001-a001*c100))*g1*h1
+(a010*b001-a001*b010)*(a100*b001-a001*b100)*(h1)^2 (5.3.15)

x*z=( b010*c001-b001*c010)*( a100*b010-a010*b100)*(f1)^2

((b010*c001-b001*c010)*(a010*c100-a100*c010+( a100*b010-a010*b100)* (-a010*c001+a001*c010))*f1*g1
+((b010*c001-b001*c010)*(a100*b010-a010*b100+(a010*b001-a001*b010)*( a100*b010-a010*b100))*f1*h1
+(-a010*c001+a001*c010)* (a010*c100-a100*c010)*(g1)^2
+((-a010*c001+a001*c010)*( a100*b010-a010*b100)+
+(a010*b001-a001*b010)* (a010*c100-a100*c010))*g1*h1
+(a010*b001-a001*b010)*(a100*b010-a010*b100)*(h1)^2 (5.3.16)

(y)^2=(-b100*c001+b001*c100)^2*(f1)^2
+2*(-b100*c001+b001*c100)*( a100*c001-a001*c100)*f1*g1
+2*(a100*b001-a001*b100)* (-b100*c001+b001*c100)*f1*h1
+(a100*c001-a001*c100)^2*(g1)^2
+2*(a100*c001-a001*c100)*( a100*b001-a001*b100)*g1*h1
+(a100*b001-a001*b100)^2*(h1)^2 (5.3.17)

y*z=(-b100*c001+b001*c100)*( b100*c010-b010*c100)*(f1)^2
+((-b100*c001+b001*c100)* (a010*c100-a100*c010)
+(a100*c001-a001*c100)*( b100*c010-b010*c100))*f1*g1
+((-b100*c001+b001*c100)* (a100*b010-a010*b100)
+(a100*b001-a001*b100)*( b100*c010-b010*c100))*f1*h1
+(a100*c001-a001*c100)* (a010*c100-a100*c010)*(g1)^2
+((a100*b010-a010*b100)* (a100*c001-a001*c100)
+(a100*b001-a001*b100)* (a010*c100-a100*c010))*g1*h1
+(a100*b001-a001*b100)*(a100*b010-a010*b100)*(h1)^2 (5.3.18)

(z)^2=(b100*c010-b010*c100)^2*(f1)^2
+2*( b100*c010-b010*c100)*( a010*c100-a100*c010)*f1*g1
+2*( b100*c010-b010*c100)* (a100*b010-a010*b100)*f1*h1
+(a010*c100-a100*c010)^2*(g1)^2
+2*(a010*c100-a100*c010)*(a100*b010-a010*b100)*g1*h1
+(a100*b010-a010*b100)^2*(h1)^2 (5.3.19)

As done previously, from these equations we can form the following single matrix equation, namely,

$$\begin{bmatrix} (x)^\wedge 2 \\ (x*y) \\ (x*z) \\ (y)^\wedge 2 \\ (y*z) \\ (z)^\wedge 2 \end{bmatrix} = B \begin{bmatrix} (f1)^\wedge 2 \\ (f1*g1) \\ (f1*h1) \\ (g1)^\wedge 2 \\ (g1*h1) \\ (h1)^\wedge 2 \end{bmatrix}$$

From the process of formation of the above matrix equations (or by actually checking) it is easy to see that $B = A^{-1}$. By proceeding for proof along the same line as is done for theorem 2.3.1 in section 2.3 the above given theorem 5.3.1 follows straightforwardly.
Instead of using symbols $f, g, h$ for polynomials and $x, y, z$ for variables if we denote three polynomials by $u_1, u_2, u_3$ and variables by $x_1, x_2, x_3$ with coefficients in the field $k$ of characteristic zero,
For $j$ = 1 to 3, we will have

$$u_j = \sum_{m_1, m_2, m_3 = 0}^{l_1(u_j), l_2(u_j), l_3(u_j)} a^j_{m_1 m_2 m_3} x_1^{m_1} x_2^{m_2} x_3^{m_3} \quad (5.3.20)$$

such that $a^j_{m_1 m_2 m_3} \in k$ and the Jacobian,

$$J_{(x_1, x_2, x_3)}(u_1, u_2, u_3) = \det \begin{pmatrix} \frac{\partial u_1}{\partial x_1} & \frac{\partial u_1}{\partial x_2} & \frac{\partial u_1}{\partial x_3} \\ \frac{\partial u_2}{\partial x_1} & \frac{\partial u_2}{\partial x_2} & \frac{\partial u_2}{\partial x_3} \\ \frac{\partial u_3}{\partial x_1} & \frac{\partial u_3}{\partial x_2} & \frac{\partial u_3}{\partial x_3} \end{pmatrix} = M = 1 \neq 0$$

By inverse function theorem we will have

$$\begin{pmatrix} \frac{\partial x_1}{\partial u_1} & \frac{\partial x_1}{\partial u_2} & \frac{\partial x_1}{\partial u_3} \\ \frac{\partial x_2}{\partial u_1} & \frac{\partial x_2}{\partial u_2} & \frac{\partial x_2}{\partial u_3} \\ \frac{\partial x_3}{\partial u_1} & \frac{\partial x_3}{\partial u_2} & \frac{\partial x_3}{\partial u_3} \end{pmatrix} = \begin{pmatrix} (\frac{\partial u_2}{\partial x_2}\frac{\partial u_3}{\partial x_3} - \frac{\partial u_2}{\partial x_3}\frac{\partial u_3}{\partial x_2}) & (\frac{\partial u_2}{\partial x_3}\frac{\partial u_3}{\partial x_1} - \frac{\partial u_2}{\partial x_1}\frac{\partial u_3}{\partial x_3}) & (\frac{\partial u_2}{\partial x_1}\frac{\partial u_3}{\partial x_2} - \frac{\partial u_2}{\partial x_2}\frac{\partial u_3}{\partial x_1}) \\ (\frac{\partial u_1}{\partial x_3}\frac{\partial u_3}{\partial x_2} - \frac{\partial u_1}{\partial x_2}\frac{\partial u_3}{\partial x_3}) & (\frac{\partial u_1}{\partial x_1}\frac{\partial u_3}{\partial x_3} - \frac{\partial u_1}{\partial x_3}\frac{\partial u_3}{\partial x_1}) & (\frac{\partial u_1}{\partial x_2}\frac{\partial u_3}{\partial x_1} - \frac{\partial u_1}{\partial x_1}\frac{\partial u_3}{\partial x_2}) \\ (\frac{\partial u_1}{\partial x_2}\frac{\partial u_2}{\partial x_3} - \frac{\partial u_1}{\partial x_3}\frac{\partial u_2}{\partial x_2}) & (\frac{\partial u_1}{\partial x_3}\frac{\partial u_2}{\partial x_1} - \frac{\partial u_1}{\partial x_1}\frac{\partial u_2}{\partial x_3}) & (\frac{\partial u_1}{\partial x_1}\frac{\partial u_2}{\partial x_2} - \frac{\partial u_1}{\partial x_2}\frac{\partial u_2}{\partial x_1}) \end{pmatrix}^T$$

(5.3.21)

where $T$ on the right hand corner indicates the operation of taking transpose and the Jacobi conditions formed as previous will be

**5.4 The Jacobi Conditions:** As for the two variables case, we can obtain on similar lines the principle as well as the derived Jacobi conditions of all types for three variables. We have the following principle Jacobi condition:

$$J_{(x_1,x_2,x_3)}(u_1,u_2,u_3) = \det\begin{pmatrix} a^1_{100} & a^1_{010} & a^1_{001} \\ a^2_{100} & a^2_{010} & a^2_{001} \\ a^3_{100} & a^3_{010} & a^3_{001} \end{pmatrix} = M = 1 \neq 0 \quad (5.4.1)$$

Collecting the coefficients of $x_1$ we have the following first derived Jacobi condition:

$$2a^1_{200}a^2_{010}a^3_{001} - a^1_{110}a^2_{100}a^3_{001} + a^1_{100}a^2_{110}a^3_{001} - 2a^1_{010}a^2_{200}a^3_{001}$$
$$- 2a^1_{200}a^2_{001}a^3_{010} + a^1_{101}a^2_{100}a^3_{010} - \cdots = 0,$$

i.e.

$$(2a^1_{200}a^2_{010} - a^1_{110}a^2_{100})a^3_{001} + (a^1_{100}a^2_{110} - 2a^1_{010}a^2_{200})a^3_{001} - \cdots = 0$$

(5.4.2)

etc.

**5.5. The nonsingular nature of the Matrices:** As is done in section 4.1 we can directly show the nonsingular nature of the matrices involved in the three variables case.

Let

$$\frac{a^1_{100}}{a^3_{100}} = u_1,\ \frac{a^2_{100}}{a^3_{100}} = u_2,\ \frac{a^1_{010}}{a^3_{010}} = v_1,\ \frac{a^2_{100}}{a^3_{010}} = v_2,\ \frac{a^1_{001}}{a^3_{001}} = w_1,\ \frac{a^2_{001}}{a^3_{001}} = w_2 .$$

With this substitution, we have

$$\det\begin{pmatrix} a^1_{100} & a^2_{100} & a^3_{100} \\ a^1_{010} & a^2_{010} & a^3_{010} \\ a^1_{001} & a^2_{001} & a^3_{001} \end{pmatrix} = \det((a^3_{100}a^3_{010}a^3_{001})\begin{pmatrix} u_1 & u_2 & 1 \\ v_1 & v_2 & 1 \\ w_1 & w_2 & 1 \end{pmatrix}) = M = 1 \neq 0.$$

Thus, the basic matrix we get for this three variables case is $\begin{pmatrix} u_1 & u_2 & 1 \\ v_1 & v_2 & 1 \\ w_1 & w_2 & 1 \end{pmatrix}$.

Using the above substitutions for the matrix $A$ above, one can check that

$$\det(A) = M^4 = 1 \neq 0 .$$

We can now prove the following

**Theorem 5.5.1:** The determinant of the matrices **like** $U_n^T$, $n = 2,3,\cdots$ that we get for the three variables case are power of Jacobian $M$, i.e. $\det(U_n^T) = M^k$ for some positive integer $k$. Moreover, the values $k$ considered in succession form a line (**fourth diagonal line**) in the Pascal triangle.

Thus, for three variables case:

| $n$ (size) | 3 | 6 | 10 | 15 | 21 | $\cdots$ |
|---|---|---|---|---|---|---|
| $k$ (power) | 1 | 4 | 10 | 20 | 35 | $\cdots$ |

□

**5.6. The System of Equations:** For this section it is convenient to proceed in the same earlier notation, i.e. we denote polynomials by $f, g, h$ and the variables by $x, y, z$. As is done in earlier sections, especially in section 3.1, we proceed now with the three variable case and show that it identically follows like the two variable case. We show how we get the desired inverse polynomials, i.e. how we get $x, y, z$ as polynomials in $f, g, h$ and actually obtain them!

We once again quote below the previously seen important matrix equation (5.2.2) which actually leads to nine equations which result by comparing the matrices element wise. Thus we have

$$\begin{bmatrix} x_f & x_g & x_h \\ y_f & y_g & y_h \\ z_f & z_g & z_h \end{bmatrix} = \begin{bmatrix} g_y h_z - g_z h_y & f_z h_y - f_y h_z & f_y g_z - f_z g_y \\ g_x h_z - g_z h_x & f_x h_z - f_z h_x & f_z g_x - f_x g_z \\ g_x h_y - g_y h_x & f_y h_x - f_x h_y & f_x g_y - f_y g_x \end{bmatrix}$$

By comparing element wise we get following nine equations:

xf=(gy*hz-gz*hy)

xg=(fz*hy-fy*hz)

xh=(fy*gz-fz*gy)

yf=(gz*hx-gx*hz)

yg=(fx*hz-fz*hx)

yh=(fz*gx-fx*gz)

zf=(gx*hy-gy*hx)

zg=(fy*hx-fx*hy)

zh=(fx*gy-fy*gx)

The **right hand side** of the above equations is **completely known** to us (given) since we are given polynomials f, g ,h in terms of variables x, y, z and we can the find the respective derivatives and can completely evaluate and determine the required expressions on the right hand side. On the other hand the **left hand side** is **completely unknown.** Here, we have to take x, y, z as power series with unknown coefficients in the basic variables f, g, h. Thus, on the left hand side we have partial derivatives of x, y, z and these x, y, z have been expressed as power series in f, g, h with unknown coefficients and our problem here is two fold: (1) We have to show that these power series are actually polynomials and (2) We have to find these unknown coefficients in terms of the known coefficients of originally given polynomials f, g, h.

**Let us first proceed to find the expressions for partial derivatives on left hand side in the above equations in terms of known right hand side:**

Let us suppose that we are given the following cubic polynomials f, g, h:

f=a100*x+a010*y+a001*z+a200*x^2+a110*x*y+a101*x*z+a020*y^2+a011*y*z+a002*z^2+a300*x^3+a210*x^2*y+a201*x^2*z+a120*x*y^2+a111*x*y*z+a102*x*z^2+a030*y^3+a021*y^2*z+a012*y*z^2+a003*z^3 (5.6.1)

g=b100*x+b010*y+b001*z+b200*x^2+b110*x*y+b101*x*z+b020*y^2+b011*y*z+b002*z^2+b300*x^3+b210*x^2*y+b201*x^2*z+b120*x*y^2+b111*x*y*z+b102*x*z^2+b030*y^3+b021*y^2*z+b012*y*z^2+b003*z^3 (5.6.2)

h=c100*x+c010*y+c001*z+c200*x^2+c110*x*y+c101*x*z+c020*y^2+c011*y*z+c002*z^2+c300*x^3+c210*x^2*y+c201*x^2*z+c120*x*y^2+c111*x*y*z+c102*x*z^2+c030*y^3+c021*y^2*z+c012*y*z^2+c003*z^3 (5.6.3)

xf=(gy*hz-gz*hy)=
(b010+b110*x+2*b020*y+b011*z+b210*x^2+2*b120*x*y+b111*x*z+3*b030*y^2+2*b021*y*z+b012*z^2)*(c001+c101*x+c011*y+2*c002*z+c201*x^2+c111*x*y+2*c102*x*z+c021*y^2+2*c012*y*z+3*c003*z^2)-
(b001+b101*x+b011*y+2*b002*z+b201*x^2+b111*x*y+2*b102

*x*z+b021*y^2+2*b012*y*z+3*b003*z^2)*(c010+c110*x+2*c0
20*y+c011*z+c210*x^2+2*c120*x*y+c111*x*z+3*c030*y^2+2*
c021*y*z+c012*z^2) (5.6.4)

xg=(fz*hy-fy*hz)=
(a001+a101*x+a011*y+2*a002*z+a201*x^2+a111*x*y+2*a102*
x*z+a021*y^2+2*a012*y*z+3*a003*z^2)*(c010+c110*x+2*c020
*y+c011*z+c210*x^2+2*c120*x*y+c111*x*z+3*c030*y^2+2*c0
21*y*z+c012*z^2)-
(a010+a110*x+2*a020*y+a011*z+a210*x^2+2*a120*x*y+a111*
x*z+3*a030*y^2+2*a021*y*z+a012*z^2)*(c001+c101*x+c011*y
+2*c002*z+c201*x^2+c111*x*y+2*c102*x*z+c021*y^2+2*c012
*y*z+3*c003*z^2) (5.6.5)

xh=(fy*gz-fz*gy)=
(a010+a110*x+2*a020*y+a011*z+a210*x^2+2*a120*x*y+a111*
x*z+3*a030*y^2+2*a021*y*z+a012*z^2)*(b001+b101*x+b011*
y+2*b002*z+b201*x^2+b111*x*y+2*b102*x*z+b021*y^2+2*b0
12*y*z+3*b003*z^2)-
(a001+a101*x+a011*y+2*a002*z+a201*x^2+a111*x*y+2*a102*
x*z+a021*y^2+2*a012*y*z+3*a003*z^2)*(b010+b110*x+2*b02
0*y+b011*z+b210*x^2+2*b120*x*y+b111*x*z+3*b030*y^2+2*
b021*y*z+b012*z^2) (5.6.6)

yf=(gz*hx-gx*hz)=
(b001+b101*x+b011*y+2*b002*z+b201*x^2+b111*x*y+2*b102
*x*z+b021*y^2+2*b012*y*z+3*b003*z^2)*(c100+2*c200*x+c1
10*y+c101*z+3*c300*x^2+2*c210*x*y+2*c201*x*z+c120*y^2+
c111*y*z+c102*z^2)-
(b100+2*b200*x+b110*y+b101*z+3*b300*x^2+2*b210*x*y+2*
b201*x*z+b120*y^2+b111*y*z+b102*z^2)*(c001+c101*x+c011
*y+2*c002*z+c201*x^2+c111*x*y+2*c102*x*z+c021*y^2+2*c0
12*y*z+3*c003*z^2) (5.6.7)

yg=(fx*hz-fz*hx)=
(a100+2*a200*x+a110*y+a101*z+3*a300*x^2+2*a210*x*y+2*a201*x*z+a120*y^2+a111*y*z+a102*z^2)*(c001+c101*x+c011*y+2*c002*z+c201*x^2+c111*x*y+2*c102*x*z+c021*y^2+2*c012*y*z+3*c003*z^2)-
(a001+a101*x+a011*y+2*a002*z+a201*x^2+a111*x*y+2*a102*x*z+a021*y^2+2*a012*y*z+3*a003*z^2)*(c100+2*c200*x+c110*y+c101*z+3*c300*x^2+2*c210*x*y+2*c201*x*z+c120*y^2+c111*y*z+c102*z^2) (5.6.8)

yh=(fz*gx-fx*gz)=
(a001+a101*x+a011*y+2*a002*z+a201*x^2+a111*x*y+2*a102*x*z+a021*y^2+2*a012*y*z+3*a003*z^2)*(b100+2*b200*x+b110*y+b101*z+3*b300*x^2+2*b210*x*y+2*b201*x*z+b120*y^2+b111*y*z+b102*z^2)-
(a100+2*a200*x+a110*y+a101*z+3*a300*x^2+2*a210*x*y+2*a201*x*z+a120*y^2+a111*y*z+a102*z^2)*(b001+b101*x+b011*y+2*b002*z+b201*x^2+b111*x*y+2*b102*x*z+b021*y^2+2*b012*y*z+3*b003*z^2) (5.6.9)

zf=(gx*hy-gy*hx)=
(b100+2*b200*x+b110*y+b101*z+3*b300*x^2+2*b210*x*y+2*b201*x*z+b120*y^2+b111*y*z+b102*z^2)*(c010+c110*x+2*c020*y+c011*z+c210*x^2+2*c120*x*y+c111*x*z+3*c030*y^2+2*c021*y*z+c012*z^2)-
(b010+b110*x+2*b020*y+b011*z+b210*x^2+2*b120*x*y+b111*x*z+3*b030*y^2+2*b021*y*z+b012*z^2)*(c100+2*c200*x+c110*y+c101*z+3*c300*x^2+2*c210*x*y+2*c201*x*z+c120*y^2+c111*y*z+c102*z^2) (5.6.10)

zg=(fy*hx-fx*hy)=
(a010+a110*x+2*a020*y+a011*z+a210*x^2+2*a120*x*y+a111*x*z+3*a030*y^2+2*a021*y*z+a012*z^2)*(c100+2*c200*x+c110*y+c101*z+3*c300*x^2+2*c210*x*y+2*c201*x*z+c120*y^2+c111*y*z+c102*z^2)-

(a100+2*a200*x+a110*y+a101*z+3*a300*x^2+2*a210*x*y+2*a201*x*z+a120*y^2+a111*y*z+a102*z^2)*(c010+c110*x+2*c020*y+c011*z+c210*x^2+2*c120*x*y+c111*x*z+3*c030*y^2+2*c021*y*z+c012*z^2) (5.6.11)

zh=(fx*gy-fy*gx)=
(a100+2*a200*x+a110*y+a101*z+3*a300*x^2+2*a210*x*y+2*a201*x*z+a120*y^2+a111*y*z+a102*z^2)*(b010+b110*x+2*b020*y+b011*z+b210*x^2+2*b120*x*y+b111*x*z+3*b030*y^2+2*b021*y*z+b012*z^2)-
(a010+a110*x+2*a020*y+a011*z+a210*x^2+2*a120*x*y+a111*x*z+3*a030*y^2+2*a021*y*z+a012*z^2)*(b100+2*b200*x+b110*y+b101*z+3*b300*x^2+2*b210*x*y+2*b201*x*z+b120*y^2+b111*y*z+b102*z^2) (5.6.12)

**Let us first proceed to find the expressions for partial derivatives on left hand side in the above equations in terms of partial derivatives of power series with unknown coefficients:**

Thus, we have following power series, viz,

x=p100*f+p010*g+p001*h+p200*f^2+p110*f*g+p101*f*h+p020*g^2+p011*g*h+p002*h^2+p300*f^3+p210*f^2*g+p201*f^2*h+p120*f*g^2+p111*f*g*h+p102*f*h^2+p030*g^3+p021*g^2*h+p012*g*h^2+p003*h^3+….. (5.6.13)

y=q100*f+q010*g+q001*h+q200*f^2+q110*f*g+q101*f*h+q020*g^2+q011*g*h+q002*h^2+q300*f^3+q210*f^2*g+q201*f^2*h+q120*f*g^2+q111*f*g*h+q102*f*h^2+q030*g^3+q021*g^2*h+q012*g*h^2+q003*h^3+…… (5.6.14)

z=r100*f+r010*g+r001*h+r200*f^2+r110*f*g+r101*f*h+r020*g^2+r011*g*h+r002*h^2+r300*f^3+r210*f^2*g+r201*f^2*h+r120*f*g^2+r111*f*g*h+r102*f*h^2+r030*g^3+r021*g^2*h+r012*g*h^2+r003*h^3+……. (5.6.14)

The partial derivatives of these power series with respect to f, g, h are

xf=p100+2*p200*f+p110*g+p101*h+3*p300*f^2+2*p210*f*g+2
*p201*f*h+p120*g^2+p111*g*h+p102*h^2+….. (5.6.15)

xg=p010+p110*f+2*p020*g+p011*h+p210*f^2+2*p120*f*g+p11
1*f*h+3*p030*g^2+2*p021*g*h+p012*h^2+….. (5.6.16)

xh=p001+p101*f+p011*g+2*p002*h+p201*f^2+p111*f*g+2*p10
2*f*h+p021*g^2+2*p012*g*h+3*p003*h^2+….. (5.6.17)

yf=q100+2*q200*f+q110*g+q101*h+3*q300*f^2+2*q210*f*g+2
*q201*f*h+q120*g^2+q111*g*h+q102*h^2+…… (5.6.18)

yg=q010+q110*f+2*q020*g+q011*h+q210*f^2+2*q120*f*g+q11
1*f*h+3*q030*g^2+2*q021*g*h+q012*h^2+…… (5.6.19)

yh=q001+q101*f+q011*g+2*q002*h+q201*f^2+q111*f*g+2*q10
2*f*h+q021*g^2+2*q012*g*h+3*q003*h^2+…... (5.6.20)

zf=r100+2*r200*f+r110*g+r101*h+3*r300*f^2+2*r210*f*g+2*r
201*f*h+r120*g^2+r111*g*h+r102*h^2+….. (5.6.21)

zg=r010+r110*f+2*r020*g+r011*h+r210*f^2+2*r120*f*g+r111*
f*h+3*r030*g^2+2*r021*g*h+r012*h^2+….. (5.6.22)

zh=r001+r101*f+r011*g+2*r002*h+r201*f^2+r111*f*g+2*r102*
f*h+r021*g^2+2*r012*g*h+3*r003*h^2+….. (5.6.23)

After back substituting the (originally given) polynomial expressions for f, g, h given in equations (5.6.1), (5.6.2), (5.6.3) respectively in these above given partial derivatives of the power series we obtain:

xf=p100+2*p200*(a100*x+a010*y+a001*z+a200*x^2+a110*x*y
+a101*x*z+a020*y^2+a011*y*z+a002*z^2+a300*x^3+a210*x^2
*y+a201*x^2*z+a120*x*y^2+a111*x*y*z+a102*x*z^2+a030*y^
3+a021*y^2*z+a012*y*z^2+a003*z^3)+p110*(b100*x+b010*y+

b001*z+b200*x^2+b110*x*y+b101*x*z+b020*y^2+b011*y*z+b
002*z^2+b300*x^3+b210*x^2*y+b201*x^2*z+b120*x*y^2+b11
1*x*y*z+b102*x*z^2+b030*y^3+b021*y^2*z+b012*y*z^2+b003
*z^3)+p101*(c100*x+c010*y+c001*z+c200*x^2+c110*x*y+c10
1*x*z+c020*y^2+c011*y*z+c002*z^2+c300*x^3+c210*x^2*y+c
201*x^2*z+c120*x*y^2+c111*x*y*z+c102*x*z^2+c030*y^3+c0
21*y^2*z+c012*y*z^2+c003*z^3)+3*p300*(a100*x+a010*y+a00
1*z+a200*x^2+a110*x*y+a101*x*z+a020*y^2+a011*y*z+a002*
z^2+a300*x^3+a210*x^2*y+a201*x^2*z+a120*x*y^2+a111*x*y
*z+a102*x*z^2+a030*y^3+a021*y^2*z+a012*y*z^2+a003*z^3)^
2+2*p210*(a100*x+a010*y+a001*z+a200*x^2+a110*x*y+a101*
x*z+a020*y^2+a011*y*z+a002*z^2+a300*x^3+a210*x^2*y+a20
1*x^2*z+a120*x*y^2+a111*x*y*z+a102*x*z^2+a030*y^3+a021
*y^2*z+a012*y*z^2+a003*z^3)*(b100*x+b010*y+b001*z+b200
*x^2+b110*x*y+b101*x*z+b020*y^2+b011*y*z+b002*z^2+b30
0*x^3+b210*x^2*y+b201*x^2*z+b120*x*y^2+b111*x*y*z+b10
2*x*z^2+b030*y^3+b021*y^2*z+b012*y*z^2+b003*z^3)+2*p20
1*(a100*x+a010*y+a001*z+a200*x^2+a110*x*y+a101*x*z+a02
0*y^2+a011*y*z+a002*z^2+a300*x^3+a210*x^2*y+a201*x^2*z
+a120*x*y^2+a111*x*y*z+a102*x*z^2+a030*y^3+a021*y^2*z+
a012*y*z^2+a003*z^3)*(c100*x+c010*y+c001*z+c200*x^2+c11
0*x*y+c101*x*z+c020*y^2+c011*y*z+c002*z^2+c300*x^3+c21
0*x^2*y+c201*x^2*z+c120*x*y^2+c111*x*y*z+c102*x*z^2+c0
30*y^3+c021*y^2*z+c012*y*z^2+c003*z^3)+p120*(b100*x+b0
10*y+b001*z+b200*x^2+b110*x*y+b101*x*z+b020*y^2+b011*
y*z+b002*z^2+b300*x^3+b210*x^2*y+b201*x^2*z+b120*x*y^2
+b111*x*y*z+b102*x*z^2+b030*y^3+b021*y^2*z+b012*y*z^2
+b003*z^3)^2+p111*(b100*x+b010*y+b001*z+b200*x^2+b110*
x*y+b101*x*z+b020*y^2+b011*y*z+b002*z^2+b300*x^3+b210
*x^2*y+b201*x^2*z+b120*x*y^2+b111*x*y*z+b102*x*z^2+b0
30*y^3+b021*y^2*z+b012*y*z^2+b003*z^3)*(c100*x+c010*y+
c001*z+c200*x^2+c110*x*y+c101*x*z+c020*y^2+c011*y*z+c0
02*z^2+c300*x^3+c210*x^2*y+c201*x^2*z+c120*x*y^2+c111*
x*y*z+c102*x*z^2+c030*y^3+c021*y^2*z+c012*y*z^2+c003*z
^3)+p102*(c100*x+c010*y+c001*z+c200*x^2+c110*x*y+c101*

x*z+c020*y^2+c011*y*z+c002*z^2+c300*x^3+c210*x^2*y+c20
1*x^2*z+c120*x*y^2+c111*x*y*z+c102*x*z^2+c030*y^3+c021
*y^2*z+c012*y*z^2+c003*z^3)^2+….. (5.6.24)

xg=p010+p110*(a100*x+a010*y+a001*z+a200*x^2+a110*x*y+a
101*x*z+a020*y^2+a011*y*z+a002*z^2+a300*x^3+a210*x^2*y
+a201*x^2*z+a120*x*y^2+a111*x*y*z+a102*x*z^2+a030*y^3+
a021*y^2*z+a012*y*z^2+a003*z^3)+2*p020*(b100*x+b010*y+
b001*z+b200*x^2+b110*x*y+b101*x*z+b020*y^2+b011*y*z+b
002*z^2+b300*x^3+b210*x^2*y+b201*x^2*z+b120*x*y^2+b11
1*x*y*z+b102*x*z^2+b030*y^3+b021*y^2*z+b012*y*z^2+b003
*z^3)+p011*(c100*x+c010*y+c001*z+c200*x^2+c110*x*y+c10
1*x*z+c020*y^2+c011*y*z+c002*z^2+c300*x^3+c210*x^2*y+c
201*x^2*z+c120*x*y^2+c111*x*y*z+c102*x*z^2+c030*y^3+c0
21*y^2*z+c012*y*z^2+c003*z^3)+p210*(a100*x+a010*y+a001*
z+a200*x^2+a110*x*y+a101*x*z+a020*y^2+a011*y*z+a002*z^
2+a300*x^3+a210*x^2*y+a201*x^2*z+a120*x*y^2+a111*x*y*z
+a102*x*z^2+a030*y^3+a021*y^2*z+a012*y*z^2+a003*z^3)^2
+2*p120*(a100*x+a010*y+a001*z+a200*x^2+a110*x*y+a101*x
*z+a020*y^2+a011*y*z+a002*z^2+a300*x^3+a210*x^2*y+a201
*x^2*z+a120*x*y^2+a111*x*y*z+a102*x*z^2+a030*y^3+a021*
y^2*z+a012*y*z^2+a003*z^3)*(b100*x+b010*y+b001*z+b200*
x^2+b110*x*y+b101*x*z+b020*y^2+b011*y*z+b002*z^2+b300
*x^3+b210*x^2*y+b201*x^2*z+b120*x*y^2+b111*x*y*z+b102
*x*z^2+b030*y^3+b021*y^2*z+b012*y*z^2+b003*z^3)+p111*(
a100*x+a010*y+a001*z+a200*x^2+a110*x*y+a101*x*z+a020*y
^2+a011*y*z+a002*z^2+a300*x^3+a210*x^2*y+a201*x^2*z+a1
20*x*y^2+a111*x*y*z+a102*x*z^2+a030*y^3+a021*y^2*z+a01
2*y*z^2+a003*z^3)*(c100*x+c010*y+c001*z+c200*x^2+c110*x
*y+c101*x*z+c020*y^2+c011*y*z+c002*z^2+c300*x^3+c210*x
^2*y+c201*x^2*z+c120*x*y^2+c111*x*y*z+c102*x*z^2+c030*
y^3+c021*y^2*z+c012*y*z^2+c003*z^3)+3*p030*(b100*x+b01
0*y+b001*z+b200*x^2+b110*x*y+b101*x*z+b020*y^2+b011*y
*z+b002*z^2+b300*x^3+b210*x^2*y+b201*x^2*z+b120*x*y^2
+b111*x*y*z+b102*x*z^2+b030*y^3+b021*y^2*z+b012*y*z^2

+b003*z^3)^2+2*p021*(b100*x+b010*y+b001*z+b200*x^2+b11
0*x*y+b101*x*z+b020*y^2+b011*y*z+b002*z^2+b300*x^3+b2
10*x^2*y+b201*x^2*z+b120*x*y^2+b111*x*y*z+b102*x*z^2+
b030*y^3+b021*y^2*z+b012*y*z^2+b003*z^3)*(c100*x+c010*
y+c001*z+c200*x^2+c110*x*y+c101*x*z+c020*y^2+c011*y*z+
c002*z^2+c300*x^3+c210*x^2*y+c201*x^2*z+c120*x*y^2+c11
1*x*y*z+c102*x*z^2+c030*y^3+c021*y^2*z+c012*y*z^2+c003
*z^3)+p012*(c100*x+c010*y+c001*z+c200*x^2+c110*x*y+c10
1*x*z+c020*y^2+c011*y*z+c002*z^2+c300*x^3+c210*x^2*y+c
201*x^2*z+c120*x*y^2+c111*x*y*z+c102*x*z^2+c030*y^3+c0
21*y^2*z+c012*y*z^2+c003*z^3)^2+….. (5.6.25)

xh=p001+p101*(a100*x+a010*y+a001*z+a200*x^2+a110*x*y+a
101*x*z+a020*y^2+a011*y*z+a002*z^2+a300*x^3+a210*x^2*y
+a201*x^2*z+a120*x*y^2+a111*x*y*z+a102*x*z^2+a030*y^3+
a021*y^2*z+a012*y*z^2+a003*z^3)+p011*(b100*x+b010*y+b0
01*z+b200*x^2+b110*x*y+b101*x*z+b020*y^2+b011*y*z+b00
2*z^2+b300*x^3+b210*x^2*y+b201*x^2*z+b120*x*y^2+b111*
x*y*z+b102*x*z^2+b030*y^3+b021*y^2*z+b012*y*z^2+b003*z
^3)+2*p002*(c100*x+c010*y+c001*z+c200*x^2+c110*x*y+c10
1*x*z+c020*y^2+c011*y*z+c002*z^2+c300*x^3+c210*x^2*y+c
201*x^2*z+c120*x*y^2+c111*x*y*z+c102*x*z^2+c030*y^3+c0
21*y^2*z+c012*y*z^2+c003*z^3)+p201*(a100*x+a010*y+a001*
z+a200*x^2+a110*x*y+a101*x*z+a020*y^2+a011*y*z+a002*z^
2+a300*x^3+a210*x^2*y+a201*x^2*z+a120*x*y^2+a111*x*y*z
+a102*x*z^2+a030*y^3+a021*y^2*z+a012*y*z^2+a003*z^3)^2
+p111*(a100*x+a010*y+a001*z+a200*x^2+a110*x*y+a101*x*z
+a020*y^2+a011*y*z+a002*z^2+a300*x^3+a210*x^2*y+a201*x
^2*z+a120*x*y^2+a111*x*y*z+a102*x*z^2+a030*y^3+a021*y^
2*z+a012*y*z^2+a003*z^3)*(b100*x+b010*y+b001*z+b200*x^
2+b110*x*y+b101*x*z+b020*y^2+b011*y*z+b002*z^2+b300*x
^3+b210*x^2*y+b201*x^2*z+b120*x*y^2+b111*x*y*z+b102*x
*z^2+b030*y^3+b021*y^2*z+b012*y*z^2+b003*z^3)+2*p102*(
a100*x+a010*y+a001*z+a200*x^2+a110*x*y+a101*x*z+a020*y
^2+a011*y*z+a002*z^2+a300*x^3+a210*x^2*y+a201*x^2*z+a1

20*x*y^2+a111*x*y*z+a102*x*z^2+a030*y^3+a021*y^2*z+a01
2*y*z^2+a003*z^3)*(c100*x+c010*y+c001*z+c200*x^2+c110*x
*y+c101*x*z+c020*y^2+c011*y*z+c002*z^2+c300*x^3+c210*x
^2*y+c201*x^2*z+c120*x*y^2+c111*x*y*z+c102*x*z^2+c030*
y^3+c021*y^2*z+c012*y*z^2+c003*z^3)+p021*(b100*x+b010*
y+b001*z+b200*x^2+b110*x*y+b101*x*z+b020*y^2+b011*y*z
+b002*z^2+b300*x^3+b210*x^2*y+b201*x^2*z+b120*x*y^2+b
111*x*y*z+b102*x*z^2+b030*y^3+b021*y^2*z+b012*y*z^2+b0
03*z^3)^2+2*p012*(b100*x+b010*y+b001*z+b200*x^2+b110*x
*y+b101*x*z+b020*y^2+b011*y*z+b002*z^2+b300*x^3+b210*
x^2*y+b201*x^2*z+b120*x*y^2+b111*x*y*z+b102*x*z^2+b03
0*y^3+b021*y^2*z+b012*y*z^2+b003*z^3)*(c100*x+c010*y+c
001*z+c200*x^2+c110*x*y+c101*x*z+c020*y^2+c011*y*z+c00
2*z^2+c300*x^3+c210*x^2*y+c201*x^2*z+c120*x*y^2+c111*x
*y*z+c102*x*z^2+c030*y^3+c021*y^2*z+c012*y*z^2+c003*z^
3)+3*p003*(c100*x+c010*y+c001*z+c200*x^2+c110*x*y+c101
*x*z+c020*y^2+c011*y*z+c002*z^2+c300*x^3+c210*x^2*y+c2
01*x^2*z+c120*x*y^2+c111*x*y*z+c102*x*z^2+c030*y^3+c02
1*y^2*z+c012*y*z^2+c003*z^3)^2+….. (5.6.26)

yf=q100+2*q200*(a100*x+a010*y+a001*z+a200*x^2+a110*x*y
+a101*x*z+a020*y^2+a011*y*z+a002*z^2+a300*x^3+a210*x^2
*y+a201*x^2*z+a120*x*y^2+a111*x*y*z+a102*x*z^2+a030*y^
3+a021*y^2*z+a012*y*z^2+a003*z^3)+q110*(b100*x+b010*y+
b001*z+b200*x^2+b110*x*y+b101*x*z+b020*y^2+b011*y*z+b
002*z^2+b300*x^3+b210*x^2*y+b201*x^2*z+b120*x*y^2+b11
1*x*y*z+b102*x*z^2+b030*y^3+b021*y^2*z+b012*y*z^2+b003
*z^3)+q101*(c100*x+c010*y+c001*z+c200*x^2+c110*x*y+c10
1*x*z+c020*y^2+c011*y*z+c002*z^2+c300*x^3+c210*x^2*y+c
201*x^2*z+c120*x*y^2+c111*x*y*z+c102*x*z^2+c030*y^3+c0
21*y^2*z+c012*y*z^2+c003*z^3)+3*q300*(a100*x+a010*y+a00
1*z+a200*x^2+a110*x*y+a101*x*z+a020*y^2+a011*y*z+a002*
z^2+a300*x^3+a210*x^2*y+a201*x^2*z+a120*x*y^2+a111*x*y
*z+a102*x*z^2+a030*y^3+a021*y^2*z+a012*y*z^2+a003*z^3)^
2+2*q210*(a100*x+a010*y+a001*z+a200*x^2+a110*x*y+a101*

x*z+a020*y^2+a011*y*z+a002*z^2+a300*x^3+a210*x^2*y+a20
1*x^2*z+a120*x*y^2+a111*x*y*z+a102*x*z^2+a030*y^3+a021
*y^2*z+a012*y*z^2+a003*z^3)*(b100*x+b010*y+b001*z+b200
*x^2+b110*x*y+b101*x*z+b020*y^2+b011*y*z+b002*z^2+b30
0*x^3+b210*x^2*y+b201*x^2*z+b120*x*y^2+b111*x*y*z+b10
2*x*z^2+b030*y^3+b021*y^2*z+b012*y*z^2+b003*z^3)+2*q20
1*(a100*x+a010*y+a001*z+a200*x^2+a110*x*y+a101*x*z+a02
0*y^2+a011*y*z+a002*z^2+a300*x^3+a210*x^2*y+a201*x^2*z
+a120*x*y^2+a111*x*y*z+a102*x*z^2+a030*y^3+a021*y^2*z+
a012*y*z^2+a003*z^3)*(c100*x+c010*y+c001*z+c200*x^2+c11
0*x*y+c101*x*z+c020*y^2+c011*y*z+c002*z^2+c300*x^3+c21
0*x^2*y+c201*x^2*z+c120*x*y^2+c111*x*y*z+c102*x*z^2+c0
30*y^3+c021*y^2*z+c012*y*z^2+c003*z^3)+q120*(b100*x+b0
10*y+b001*z+b200*x^2+b110*x*y+b101*x*z+b020*y^2+b011*
y*z+b002*z^2+b300*x^3+b210*x^2*y+b201*x^2*z+b120*x*y^2
+b111*x*y*z+b102*x*z^2+b030*y^3+b021*y^2*z+b012*y*z^2
+b003*z^3)^2+q111*(b100*x+b010*y+b001*z+b200*x^2+b110*
x*y+b101*x*z+b020*y^2+b011*y*z+b002*z^2+b300*x^3+b210
*x^2*y+b201*x^2*z+b120*x*y^2+b111*x*y*z+b102*x*z^2+b0
30*y^3+b021*y^2*z+b012*y*z^2+b003*z^3)*(c100*x+c010*y+
c001*z+c200*x^2+c110*x*y+c101*x*z+c020*y^2+c011*y*z+c0
02*z^2+c300*x^3+c210*x^2*y+c201*x^2*z+c120*x*y^2+c111*
x*y*z+c102*x*z^2+c030*y^3+c021*y^2*z+c012*y*z^2+c003*z
^3)+q102*(c100*x+c010*y+c001*z+c200*x^2+c110*x*y+c101*
x*z+c020*y^2+c011*y*z+c002*z^2+c300*x^3+c210*x^2*y+c20
1*x^2*z+c120*x*y^2+c111*x*y*z+c102*x*z^2+c030*y^3+c021
*y^2*z+c012*y*z^2+c003*z^3)^2+….. (5.6.27)

yg=q010+q110*(a100*x+a010*y+a001*z+a200*x^2+a110*x*y+a
101*x*z+a020*y^2+a011*y*z+a002*z^2+a300*x^3+a210*x^2*y
+a201*x^2*z+a120*x*y^2+a111*x*y*z+a102*x*z^2+a030*y^3+
a021*y^2*z+a012*y*z^2+a003*z^3)+2*q020*(b100*x+b010*y+
b001*z+b200*x^2+b110*x*y+b101*x*z+b020*y^2+b011*y*z+b
002*z^2+b300*x^3+b210*x^2*y+b201*x^2*z+b120*x*y^2+b11
1*x*y*z+b102*x*z^2+b030*y^3+b021*y^2*z+b012*y*z^2+b003

*z^3)+q011*(c100*x+c010*y+c001*z+c200*x^2+c110*x*y+c10
1*x*z+c020*y^2+c011*y*z+c002*z^2+c300*x^3+c210*x^2*y+c
201*x^2*z+c120*x*y^2+c111*x*y*z+c102*x*z^2+c030*y^3+c0
21*y^2*z+c012*y*z^2+c003*z^3)+q210*(a100*x+a010*y+a001*
z+a200*x^2+a110*x*y+a101*x*z+a020*y^2+a011*y*z+a002*z^
2+a300*x^3+a210*x^2*y+a201*x^2*z+a120*x*y^2+a111*x*y*z
+a102*x*z^2+a030*y^3+a021*y^2*z+a012*y*z^2+a003*z^3)^2
+2*q120*(a100*x+a010*y+a001*z+a200*x^2+a110*x*y+a101*x
*z+a020*y^2+a011*y*z+a002*z^2+a300*x^3+a210*x^2*y+a201
*x^2*z+a120*x*y^2+a111*x*y*z+a102*x*z^2+a030*y^3+a021*
y^2*z+a012*y*z^2+a003*z^3)*(b100*x+b010*y+b001*z+b200*
x^2+b110*x*y+b101*x*z+b020*y^2+b011*y*z+b002*z^2+b300
*x^3+b210*x^2*y+b201*x^2*z+b120*x*y^2+b111*x*y*z+b102
*x*z^2+b030*y^3+b021*y^2*z+b012*y*z^2+b003*z^3)+q111*(
a100*x+a010*y+a001*z+a200*x^2+a110*x*y+a101*x*z+a020*y
^2+a011*y*z+a002*z^2+a300*x^3+a210*x^2*y+a201*x^2*z+a1
20*x*y^2+a111*x*y*z+a102*x*z^2+a030*y^3+a021*y^2*z+a01
2*y*z^2+a003*z^3)*(c100*x+c010*y+c001*z+c200*x^2+c110*x
*y+c101*x*z+c020*y^2+c011*y*z+c002*z^2+c300*x^3+c210*x
^2*y+c201*x^2*z+c120*x*y^2+c111*x*y*z+c102*x*z^2+c030*
y^3+c021*y^2*z+c012*y*z^2+c003*z^3)+3*q030*(b100*x+b01
0*y+b001*z+b200*x^2+b110*x*y+b101*x*z+b020*y^2+b011*y
*z+b002*z^2+b300*x^3+b210*x^2*y+b201*x^2*z+b120*x*y^2
+b111*x*y*z+b102*x*z^2+b030*y^3+b021*y^2*z+b012*y*z^2
+b003*z^3)^2+2*q021*(b100*x+b010*y+b001*z+b200*x^2+b11
0*x*y+b101*x*z+b020*y^2+b011*y*z+b002*z^2+b300*x^3+b2
10*x^2*y+b201*x^2*z+b120*x*y^2+b111*x*y*z+b102*x*z^2+
b030*y^3+b021*y^2*z+b012*y*z^2+b003*z^3)*(c100*x+c010*
y+c001*z+c200*x^2+c110*x*y+c101*x*z+c020*y^2+c011*y*z+
c002*z^2+c300*x^3+c210*x^2*y+c201*x^2*z+c120*x*y^2+c11
1*x*y*z+c102*x*z^2+c030*y^3+c021*y^2*z+c012*y*z^2+c003
*z^3)+q012*(c100*x+c010*y+c001*z+c200*x^2+c110*x*y+c10
1*x*z+c020*y^2+c011*y*z+c002*z^2+c300*x^3+c210*x^2*y+c
201*x^2*z+c120*x*y^2+c111*x*y*z+c102*x*z^2+c030*y^3+c0
21*y^2*z+c012*y*z^2+c003*z^3)^2+….. (5.6.28)

yh=q001+q101*(a100*x+a010*y+a001*z+a200*x^2+a110*x*y+a
101*x*z+a020*y^2+a011*y*z+a002*z^2+a300*x^3+a210*x^2*y
+a201*x^2*z+a120*x*y^2+a111*x*y*z+a102*x*z^2+a030*y^3+
a021*y^2*z+a012*y*z^2+a003*z^3)+q011*(b100*x+b010*y+b0
01*z+b200*x^2+b110*x*y+b101*x*z+b020*y^2+b011*y*z+b00
2*z^2+b300*x^3+b210*x^2*y+b201*x^2*z+b120*x*y^2+b111*
x*y*z+b102*x*z^2+b030*y^3+b021*y^2*z+b012*y*z^2+b003*z
^3)+2*q002*(c100*x+c010*y+c001*z+c200*x^2+c110*x*y+c10
1*x*z+c020*y^2+c011*y*z+c002*z^2+c300*x^3+c210*x^2*y+c
201*x^2*z+c120*x*y^2+c111*x*y*z+c102*x*z^2+c030*y^3+c0
21*y^2*z+c012*y*z^2+c003*z^3)+q201*(a100*x+a010*y+a001*
z+a200*x^2+a110*x*y+a101*x*z+a020*y^2+a011*y*z+a002*z^
2+a300*x^3+a210*x^2*y+a201*x^2*z+a120*x*y^2+a111*x*y*z
+a102*x*z^2+a030*y^3+a021*y^2*z+a012*y*z^2+a003*z^3)^2
+q111*(a100*x+a010*y+a001*z+a200*x^2+a110*x*y+a101*x*z
+a020*y^2+a011*y*z+a002*z^2+a300*x^3+a210*x^2*y+a201*x
^2*z+a120*x*y^2+a111*x*y*z+a102*x*z^2+a030*y^3+a021*y^
2*z+a012*y*z^2+a003*z^3)*(b100*x+b010*y+b001*z+b200*x^
2+b110*x*y+b101*x*z+b020*y^2+b011*y*z+b002*z^2+b300*x
^3+b210*x^2*y+b201*x^2*z+b120*x*y^2+b111*x*y*z+b102*x
*z^2+b030*y^3+b021*y^2*z+b012*y*z^2+b003*z^3)+2*q102*(
a100*x+a010*y+a001*z+a200*x^2+a110*x*y+a101*x*z+a020*y
^2+a011*y*z+a002*z^2+a300*x^3+a210*x^2*y+a201*x^2*z+a1
20*x*y^2+a111*x*y*z+a102*x*z^2+a030*y^3+a021*y^2*z+a01
2*y*z^2+a003*z^3)*(c100*x+c010*y+c001*z+c200*x^2+c110*x
*y+c101*x*z+c020*y^2+c011*y*z+c002*z^2+c300*x^3+c210*x
^2*y+c201*x^2*z+c120*x*y^2+c111*x*y*z+c102*x*z^2+c030*
y^3+c021*y^2*z+c012*y*z^2+c003*z^3)+q021*(b100*x+b010*
y+b001*z+b200*x^2+b110*x*y+b101*x*z+b020*y^2+b011*y*z
+b002*z^2+b300*x^3+b210*x^2*y+b201*x^2*z+b120*x*y^2+b
111*x*y*z+b102*x*z^2+b030*y^3+b021*y^2*z+b012*y*z^2+b0
03*z^3)^2+2*q012*(b100*x+b010*y+b001*z+b200*x^2+b110*x
*y+b101*x*z+b020*y^2+b011*y*z+b002*z^2+b300*x^3+b210*
x^2*y+b201*x^2*z+b120*x*y^2+b111*x*y*z+b102*x*z^2+b03

0*y^3+b021*y^2*z+b012*y*z^2+b003*z^3)*(c100*x+c010*y+c
001*z+c200*x^2+c110*x*y+c101*x*z+c020*y^2+c011*y*z+c00
2*z^2+c300*x^3+c210*x^2*y+c201*x^2*z+c120*x*y^2+c111*x
*y*z+c102*x*z^2+c030*y^3+c021*y^2*z+c012*y*z^2+c003*z^
3)+3*q003*(c100*x+c010*y+c001*z+c200*x^2+c110*x*y+c101
*x*z+c020*y^2+c011*y*z+c002*z^2+c300*x^3+c210*x^2*y+c2
01*x^2*z+c120*x*y^2+c111*x*y*z+c102*x*z^2+c030*y^3+c02
1*y^2*z+c012*y*z^2+c003*z^3)^2+….. (5.6.29)

zf=r100+2*r200*(a100*x+a010*y+a001*z+a200*x^2+a110*x*y+
a101*x*z+a020*y^2+a011*y*z+a002*z^2+a300*x^3+a210*x^2*
y+a201*x^2*z+a120*x*y^2+a111*x*y*z+a102*x*z^2+a030*y^3
+a021*y^2*z+a012*y*z^2+a003*z^3)+r110*(b100*x+b010*y+b0
01*z+b200*x^2+b110*x*y+b101*x*z+b020*y^2+b011*y*z+b00
2*z^2+b300*x^3+b210*x^2*y+b201*x^2*z+b120*x*y^2+b111*
x*y*z+b102*x*z^2+b030*y^3+b021*y^2*z+b012*y*z^2+b003*z
^3)+r101*(c100*x+c010*y+c001*z+c200*x^2+c110*x*y+c101*x
*z+c020*y^2+c011*y*z+c002*z^2+c300*x^3+c210*x^2*y+c201
*x^2*z+c120*x*y^2+c111*x*y*z+c102*x*z^2+c030*y^3+c021*
y^2*z+c012*y*z^2+c003*z^3)+3*r300*(a100*x+a010*y+a001*z
+a200*x^2+a110*x*y+a101*x*z+a020*y^2+a011*y*z+a002*z^2
+a300*x^3+a210*x^2*y+a201*x^2*z+a120*x*y^2+a111*x*y*z+
a102*x*z^2+a030*y^3+a021*y^2*z+a012*y*z^2+a003*z^3)^2+2
*r210*(a100*x+a010*y+a001*z+a200*x^2+a110*x*y+a101*x*z
+a020*y^2+a011*y*z+a002*z^2+a300*x^3+a210*x^2*y+a201*x
^2*z+a120*x*y^2+a111*x*y*z+a102*x*z^2+a030*y^3+a021*y^
2*z+a012*y*z^2+a003*z^3)*(b100*x+b010*y+b001*z+b200*x^
2+b110*x*y+b101*x*z+b020*y^2+b011*y*z+b002*z^2+b300*x
^3+b210*x^2*y+b201*x^2*z+b120*x*y^2+b111*x*y*z+b102*x
*z^2+b030*y^3+b021*y^2*z+b012*y*z^2+b003*z^3)+2*r201*(a
100*x+a010*y+a001*z+a200*x^2+a110*x*y+a101*x*z+a020*y^
2+a011*y*z+a002*z^2+a300*x^3+a210*x^2*y+a201*x^2*z+a12
0*x*y^2+a111*x*y*z+a102*x*z^2+a030*y^3+a021*y^2*z+a012
*y*z^2+a003*z^3)*(c100*x+c010*y+c001*z+c200*x^2+c110*x*
y+c101*x*z+c020*y^2+c011*y*z+c002*z^2+c300*x^3+c210*x^

2*y+c201*x^2*z+c120*x*y^2+c111*x*y*z+c102*x*z^2+c030*y
^3+c021*y^2*z+c012*y*z^2+c003*z^3)+r120*(b100*x+b010*y+
b001*z+b200*x^2+b110*x*y+b101*x*z+b020*y^2+b011*y*z+b
002*z^2+b300*x^3+b210*x^2*y+b201*x^2*z+b120*x*y^2+b11
1*x*y*z+b102*x*z^2+b030*y^3+b021*y^2*z+b012*y*z^2+b003
*z^3)^2+r111*(b100*x+b010*y+b001*z+b200*x^2+b110*x*y+b
101*x*z+b020*y^2+b011*y*z+b002*z^2+b300*x^3+b210*x^2*
y+b201*x^2*z+b120*x*y^2+b111*x*y*z+b102*x*z^2+b030*y^
3+b021*y^2*z+b012*y*z^2+b003*z^3)*(c100*x+c010*y+c001*
z+c200*x^2+c110*x*y+c101*x*z+c020*y^2+c011*y*z+c002*z^
2+c300*x^3+c210*x^2*y+c201*x^2*z+c120*x*y^2+c111*x*y*z
+c102*x*z^2+c030*y^3+c021*y^2*z+c012*y*z^2+c003*z^3)+r1
02*(c100*x+c010*y+c001*z+c200*x^2+c110*x*y+c101*x*z+c0
20*y^2+c011*y*z+c002*z^2+c300*x^3+c210*x^2*y+c201*x^2*
z+c120*x*y^2+c111*x*y*z+c102*x*z^2+c030*y^3+c021*y^2*z
+c012*y*z^2+c003*z^3)^2+….. (5.6.30)

zg=r010+r110*(a100*x+a010*y+a001*z+a200*x^2+a110*x*y+a1
01*x*z+a020*y^2+a011*y*z+a002*z^2+a300*x^3+a210*x^2*y+
a201*x^2*z+a120*x*y^2+a111*x*y*z+a102*x*z^2+a030*y^3+a
021*y^2*z+a012*y*z^2+a003*z^3)+2*r020*(b100*x+b010*y+b0
01*z+b200*x^2+b110*x*y+b101*x*z+b020*y^2+b011*y*z+b00
2*z^2+b300*x^3+b210*x^2*y+b201*x^2*z+b120*x*y^2+b111*
x*y*z+b102*x*z^2+b030*y^3+b021*y^2*z+b012*y*z^2+b003*z
^3)+r011*(c100*x+c010*y+c001*z+c200*x^2+c110*x*y+c101*x
*z+c020*y^2+c011*y*z+c002*z^2+c300*x^3+c210*x^2*y+c201
*x^2*z+c120*x*y^2+c111*x*y*z+c102*x*z^2+c030*y^3+c021*
y^2*z+c012*y*z^2+c003*z^3)+r210*(a100*x+a010*y+a001*z+a
200*x^2+a110*x*y+a101*x*z+a020*y^2+a011*y*z+a002*z^2+a
300*x^3+a210*x^2*y+a201*x^2*z+a120*x*y^2+a111*x*y*z+a1
02*x*z^2+a030*y^3+a021*y^2*z+a012*y*z^2+a003*z^3)^2+2*r
120*(a100*x+a010*y+a001*z+a200*x^2+a110*x*y+a101*x*z+a
020*y^2+a011*y*z+a002*z^2+a300*x^3+a210*x^2*y+a201*x^2
*z+a120*x*y^2+a111*x*y*z+a102*x*z^2+a030*y^3+a021*y^2*
z+a012*y*z^2+a003*z^3)*(b100*x+b010*y+b001*z+b200*x^2+

b110*x*y+b101*x*z+b020*y^2+b011*y*z+b002*z^2+b300*x^3
+b210*x^2*y+b201*x^2*z+b120*x*y^2+b111*x*y*z+b102*x*z
^2+b030*y^3+b021*y^2*z+b012*y*z^2+b003*z^3)+r111*(a100*
x+a010*y+a001*z+a200*x^2+a110*x*y+a101*x*z+a020*y^2+a0
11*y*z+a002*z^2+a300*x^3+a210*x^2*y+a201*x^2*z+a120*x*
y^2+a111*x*y*z+a102*x*z^2+a030*y^3+a021*y^2*z+a012*y*z
^2+a003*z^3)*(c100*x+c010*y+c001*z+c200*x^2+c110*x*y+c1
01*x*z+c020*y^2+c011*y*z+c002*z^2+c300*x^3+c210*x^2*y+
c201*x^2*z+c120*x*y^2+c111*x*y*z+c102*x*z^2+c030*y^3+c
021*y^2*z+c012*y*z^2+c003*z^3)+3*r030*(b100*x+b010*y+b0
01*z+b200*x^2+b110*x*y+b101*x*z+b020*y^2+b011*y*z+b00
2*z^2+b300*x^3+b210*x^2*y+b201*x^2*z+b120*x*y^2+b111*
x*y*z+b102*x*z^2+b030*y^3+b021*y^2*z+b012*y*z^2+b003*z
^3)^2+2*r021*(b100*x+b010*y+b001*z+b200*x^2+b110*x*y+b
101*x*z+b020*y^2+b011*y*z+b002*z^2+b300*x^3+b210*x^2*
y+b201*x^2*z+b120*x*y^2+b111*x*y*z+b102*x*z^2+b030*y^
3+b021*y^2*z+b012*y*z^2+b003*z^3)*(c100*x+c010*y+c001*
z+c200*x^2+c110*x*y+c101*x*z+c020*y^2+c011*y*z+c002*z^
2+c300*x^3+c210*x^2*y+c201*x^2*z+c120*x*y^2+c111*x*y*z
+c102*x*z^2+c030*y^3+c021*y^2*z+c012*y*z^2+c003*z^3)+r0
12*(c100*x+c010*y+c001*z+c200*x^2+c110*x*y+c101*x*z+c0
20*y^2+c011*y*z+c002*z^2+c300*x^3+c210*x^2*y+c201*x^2*
z+c120*x*y^2+c111*x*y*z+c102*x*z^2+c030*y^3+c021*y^2*z
+c012*y*z^2+c003*z^3)^2+….. (5.6.31)

zh=r001+r101*(a100*x+a010*y+a001*z+a200*x^2+a110*x*y+a1
01*x*z+a020*y^2+a011*y*z+a002*z^2+a300*x^3+a210*x^2*y+
a201*x^2*z+a120*x*y^2+a111*x*y*z+a102*x*z^2+a030*y^3+a
021*y^2*z+a012*y*z^2+a003*z^3)+r011*(b100*x+b010*y+b001
*z+b200*x^2+b110*x*y+b101*x*z+b020*y^2+b011*y*z+b002*
z^2+b300*x^3+b210*x^2*y+b201*x^2*z+b120*x*y^2+b111*x*
y*z+b102*x*z^2+b030*y^3+b021*y^2*z+b012*y*z^2+b003*z^3
)+2*r002*(c100*x+c010*y+c001*z+c200*x^2+c110*x*y+c101*x
*z+c020*y^2+c011*y*z+c002*z^2+c300*x^3+c210*x^2*y+c201
*x^2*z+c120*x*y^2+c111*x*y*z+c102*x*z^2+c030*y^3+c021*

y^2*z+c012*y*z^2+c003*z^3)+r201*(a100*x+a010*y+a001*z+a
200*x^2+a110*x*y+a101*x*z+a020*y^2+a011*y*z+a002*z^2+a
300*x^3+a210*x^2*y+a201*x^2*z+a120*x*y^2+a111*x*y*z+a1
02*x*z^2+a030*y^3+a021*y^2*z+a012*y*z^2+a003*z^3)^2+r11
1*(a100*x+a010*y+a001*z+a200*x^2+a110*x*y+a101*x*z+a02
0*y^2+a011*y*z+a002*z^2+a300*x^3+a210*x^2*y+a201*x^2*z
+a120*x*y^2+a111*x*y*z+a102*x*z^2+a030*y^3+a021*y^2*z+
a012*y*z^2+a003*z^3)*(b100*x+b010*y+b001*z+b200*x^2+b1
10*x*y+b101*x*z+b020*y^2+b011*y*z+b002*z^2+b300*x^3+b
210*x^2*y+b201*x^2*z+b120*x*y^2+b111*x*y*z+b102*x*z^2
+b030*y^3+b021*y^2*z+b012*y*z^2+b003*z^3)+2*r102*(a100*
x+a010*y+a001*z+a200*x^2+a110*x*y+a101*x*z+a020*y^2+a0
11*y*z+a002*z^2+a300*x^3+a210*x^2*y+a201*x^2*z+a120*x*
y^2+a111*x*y*z+a102*x*z^2+a030*y^3+a021*y^2*z+a012*y*z
^2+a003*z^3)*(c100*x+c010*y+c001*z+c200*x^2+c110*x*y+c1
01*x*z+c020*y^2+c011*y*z+c002*z^2+c300*x^3+c210*x^2*y+
c201*x^2*z+c120*x*y^2+c111*x*y*z+c102*x*z^2+c030*y^3+c
021*y^2*z+c012*y*z^2+c003*z^3)+r021*(b100*x+b010*y+b001
*z+b200*x^2+b110*x*y+b101*x*z+b020*y^2+b011*y*z+b002*
z^2+b300*x^3+b210*x^2*y+b201*x^2*z+b120*x*y^2+b111*x*
y*z+b102*x*z^2+b030*y^3+b021*y^2*z+b012*y*z^2+b003*z^3
)^2+2*r012*(b100*x+b010*y+b001*z+b200*x^2+b110*x*y+b10
1*x*z+b020*y^2+b011*y*z+b002*z^2+b300*x^3+b210*x^2*y+
b201*x^2*z+b120*x*y^2+b111*x*y*z+b102*x*z^2+b030*y^3+
b021*y^2*z+b012*y*z^2+b003*z^3)*(c100*x+c010*y+c001*z+c
200*x^2+c110*x*y+c101*x*z+c020*y^2+c011*y*z+c002*z^2+c
300*x^3+c210*x^2*y+c201*x^2*z+c120*x*y^2+c111*x*y*z+c1
02*x*z^2+c030*y^3+c021*y^2*z+c012*y*z^2+c003*z^3)+3*r00
3*(c100*x+c010*y+c001*z+c200*x^2+c110*x*y+c101*x*z+c02
0*y^2+c011*y*z+c002*z^2+c300*x^3+c210*x^2*y+c201*x^2*z
+c120*x*y^2+c111*x*y*z+c102*x*z^2+c030*y^3+c021*y^2*z+
c012*y*z^2+c003*z^3)^2+….. (5.6.32)

We now proceed to form the actual systems of equations by comparing the coefficients of same degree, starting from degree zero and moving onwards to consider higher degree case in succession, on both sides for homogeneous blocks corresponding to identical degree using two types of expressions for partial derivatives xf, xg, xh, yf, yg, yh, zf, zg, zh. **Then by equating and comparing right hand sides of xf, xg, xh, yf, yg, yh, zf, zg, zh from equation sets (5.6.4) to (5.6.12) and from equations sets (5.6.24) to (5.6.32),** we have

(1) By comparing terms of **degree = 0,**

p100 = (b010*c001-b001*c010)

p010 = (a001*c010-a010*c001)

p001 = (a010*b001- a001*b010)

q100 = (b001*c100- b100*c001)

q010 = (a100*c001- a001*c100)

q001 = (a001*b100- a100*b001)

r100 = (b100*c010- b010*c100)

r010 = (a010*c100- a100*c010)

r001 = (a100*b010- a010*b100)

Note that in the two variables case the expressions on the right hand side in the above equations giving the values of unknown coefficients were determinantal monomials of the determinant of size one, but now in the present (three variables) case they are determinantal monomials of the determinant of size two. In the several, say (n+1), variable case they will be determinantal monomials of the determinant of size n as can be easily seen.

(2) By comparing terms of **degree = 1,** i.e. by collecting and equating the coefficients of x, y, and z in both expressions for “xf” we arrive at the following matrix equation:

$$U^T \begin{bmatrix} 2*p200 \\ p110 \\ p101 \end{bmatrix} = \begin{bmatrix} l_1^1 \\ m_1^1 \\ n_1^1 \end{bmatrix}$$

where $U^T$ is **transpose** of matrix $U$ given in section (5.3). Proceeding on similar lines and by collecting and equating the coefficients of x, y, and z in both expressions for “xg” we arrive at the following matrix equation:

$$U^T \begin{bmatrix} p110 \\ 2*p020 \\ p011 \end{bmatrix} = \begin{bmatrix} l_2^1 \\ m_2^1 \\ n_2^1 \end{bmatrix}$$

Proceeding on similar lines and by collecting and equating the coefficients of x, y, and z in both expressions for “xh” we arrive at the following matrix equation:

$$U^T \begin{bmatrix} p101 \\ p011 \\ 2*p002 \end{bmatrix} = \begin{bmatrix} l_3^1 \\ m_3^1 \\ n_3^1 \end{bmatrix}$$

Proceeding on similar lines and by collecting and equating the coefficients of x, y, and z in both expressions for “yf” we arrive at the following matrix equation:

$$U^T \begin{bmatrix} 2*q200 \\ q110 \\ q101 \end{bmatrix} = \begin{bmatrix} l_1^2 \\ m_1^2 \\ n_1^2 \end{bmatrix}$$

Proceeding on similar lines and by collecting and equating the coefficients of x, y, and z in both expressions for “yg” we arrive at the following matrix equation:

$$U^T \begin{bmatrix} q110 \\ 2*q020 \\ q011 \end{bmatrix} = \begin{bmatrix} l_2^2 \\ m_2^2 \\ n_2^2 \end{bmatrix}$$

Proceeding on similar lines and by collecting and equating the coefficients of x, y, and z in both expressions for “yh” we arrive at the following matrix equation:

$$U^T \begin{bmatrix} q101 \\ q011 \\ 2*q002 \end{bmatrix} = \begin{bmatrix} l_3^2 \\ m_3^2 \\ n_3^2 \end{bmatrix}$$

Proceeding on similar lines and by collecting and equating the coefficients of x, y, and z in both expressions for “zf” we arrive at the following matrix equation:

$$U^T \begin{bmatrix} 2*r200 \\ r110 \\ r101 \end{bmatrix} = \begin{bmatrix} l_1^3 \\ m_1^3 \\ n_1^3 \end{bmatrix}$$

Proceeding on similar lines and by collecting and equating the coefficients of x, y, and z in both expressions for “zg” we arrive at the following matrix equation:

$$U^T \begin{bmatrix} r110 \\ 2*r020 \\ r011 \end{bmatrix} = \begin{bmatrix} l_2^3 \\ m_2^3 \\ n_2^3 \end{bmatrix}$$

Proceeding on similar lines and by collecting and equating the coefficients of x, y, and z in both expressions for "zh" we arrive at the following matrix equation:

$$U^T \begin{bmatrix} r101 \\ r011 \\ 2*r002 \end{bmatrix} = \begin{bmatrix} l_3^3 \\ m_3^3 \\ n_3^3 \end{bmatrix}$$

In all these matrix equations the right hand side is completely known, i.e. all the constants $l_j^i, m_j^i, n_j^i$ for all $i, j$ are known as they all are made up from known constants, namely, the coefficients of given polynomials $f, g, h$. As an illustration we quote below the values of first few of these constants,

$$l_1^1 = (b010*c101 + b110*c001) - (b001*c110 + b101*c010)$$

$$m_2^2 = (a100*c011 + a110*c001) - (a001*c110 + a011*c100)$$

$$n_3^3 = (a100*b011 + a101*b010) - (a010*b101 + a011*b100)$$

$$\vdots$$

etc. Now, since $U$ is invertible and $U^{-1} = V$, therefore, $U^T$ is also invertible and $(U^T)^{-1} = V^T$, where $V^T$ is the transpose of matrix $V$ stated just above the statement of theorem 5.3.1. Thus, we can uniquely determine the unknown coefficients of quadratic terms from the above matrix equations.

(3) By comparing the terms of **degree = 2**, i.e. by collecting and equating the coefficients of x^2, x*y, x*z, y^2, y*z, z^2 in both expressions for "xf" we arrive at the six equations of the following type. Note that the equation given below is obtained by collecting and equating the **coefficients of x^2 in both expressions for "xf".** The term inside the square bracket, i.e. inside [ ], is called **residuum** in which each term contains a coefficient corresponding to some quadratic term, like, a200, b200, c200.

3*p300*a100^2+2*p210*a100*b100+2*p201*a100*c100
+p120*b100^2+p111*b100*c100+p102*c100^2
=((b010*c201+b110*c101+b210*c001)
-(b001*c210+b101*c110+b201*c010))
-[2*p200*a200+p110*b200+p101*c200]

Note that the equation given below is obtained by collecting and equating the **coefficients of x*y in both expressions for "xf".** The term inside the square bracket, i.e. inside [ ], is called **residuum** in which each term contains a coefficient corresponding to some quadratic term, like, a110, b110, c110.

3*p300*2*a100*a010+2*p210*(a100*b010+a010*b100)
+2*p201*(a100*c010+a010*c100)+p120*2*b100*b010
+p111*(b100*c010+b010*c100)+p102*2*c100*c010)
=((b010*c111+b110*c011+2*b120*c001+2*b020*c101)
-(b001*2*c120+b101*2*c020+b011*c110+b111*c010))
-[2*p200*a110+p110*b110+p101*c110]

We can proceed on these lines and obtain the remaining four similar equations by collecting and equating the coefficients of x*z, y^2, y*z, z^2 in both expressions for **"xf".** We can then finally collect these six equations into a single matrix equation of following type:

.

$$A^T \begin{bmatrix} 3*p300 \\ 2*p210 \\ 2*p201 \\ p120 \\ p111 \\ p102 \end{bmatrix} = \begin{bmatrix} K(xf)-Y(xf,x\wedge 2) \\ L(xf)-Y(xf,x*y) \\ M(xf)-Y(xf,x*z) \\ N(xf)-Y(xf,y\wedge 2) \\ O(xf)-Y(xf,y*z) \\ P(xf)-Y(xf,z\wedge 2) \end{bmatrix}$$

where $A^T$ is the transpose of the matrix $A$ defined in the matrix equation (5.3.13). Now, since $A$ is invertible, and $A^{-1} = B$, where $B$ is the matrix defined in the matrix equation that's given just below equation (5.3.19), therefore, $A^T$ is also invertible, $\det(A) = 1$, and $(A^T)^{-1} = B^T$, where $B^T$ is the transpose of matrix $B$. Note that quantities $K^s, L^s, M^s, N^s, O^s, P^s$ are in terms of known coefficients of $f, g, h$ and so completely known. The terms $Y^s$ are residuum terms and we will see that we can suppose that these terms are actually nonexistent (i.e. they vanish when we choose BCW form for the given cubic polynomials in which the quadratic terms are absent and a quadratic multiplier coefficient exists in every term of which $Y^s$ are made up of, therefore, we can assume that all $Y^s$ are individually zero). Also, it is **important** to note that **when given polynomial have BCW form** then quadratic terms are absent in the given polynomials, therefore, all coefficients of quadratic terms are equal to zero, i.e. all $a_{ijk} = 0, b_{ijk} = 0, c_{ijk}$ whenever $i + j + k = 2$ and this further implies from matrix equations obtained in (2) above by comparing terms of **degree = 1,** that also $p_{ijk} = 0, q_{ijk} = 0, r_{ijk} = 0$ whenever $i + j + k = 2$ **when given polynomials have BCW form.** Thus, by multiplying both sides of the above matrix equation by $(A^T)^{-1} = B^T$, i.e. by inverting the above matrix equation, we can uniquely determine the unknown cubic coefficients present there.

By comparing the terms of **degree = 2**, i.e. by collecting and equating the coefficients of x^2, x*y, x*z, y^2, y*z, z^2 in both expressions for **"xg"** we arrive at the six equations and we can collect these equations into a single matrix equation given below:

$$
A^T \begin{bmatrix} p210 \\ 2*p120 \\ p111 \\ 3*p030 \\ 2*p021 \\ p012 \end{bmatrix} = \begin{bmatrix} K(xg) - Y(xg, x\hat{}2) \\ L(xg) - Y(xg, x*y) \\ M(xg) - Y(xg, x*z) \\ N(xg) - Y(xg, y\hat{}2) \\ O(xg) - Y(xg, y*z) \\ P(xg) - Y(xg, z\hat{}2) \end{bmatrix}
$$

As previous, we can uniquely determine the unknown cubic coefficients present in the above matrix equation.
By comparing the terms of **degree = 2**, i.e. by collecting and equating the coefficients of x^2, x*y, x*z, y^2, y*z, z^2 in both expressions for **“xh”** we arrive at the six equations and we can collect these equations into a single matrix equation given below:

$$
A^T \begin{bmatrix} p201 \\ p111 \\ 2*p102 \\ p021 \\ 2*p012 \\ 3*p003 \end{bmatrix} = \begin{bmatrix} K(xh) - Y(xh, x\hat{}2) \\ L(xh) - Y(xh, x*y) \\ M(xh) - Y(xh, x*z) \\ N(xh) - Y(xh, y\hat{}2) \\ O(xh) - Y(xh, y*z) \\ P(xh) - Y(xh, z\hat{}2) \end{bmatrix}
$$

As previous, we can uniquely determine the unknown cubic coefficients present in the above matrix equation.

By comparing the terms of **degree = 2**, i.e. by collecting and equating the coefficients of x^2, x*y, x*z, y^2, y*z, z^2 in both expressions for **“yf”** we arrive at the six equations and we can collect these equations into a single matrix equation given below:

$$A^T \begin{bmatrix} 3*q300 \\ 2*q210 \\ 2*q201 \\ q120 \\ q111 \\ q102 \end{bmatrix} = \begin{bmatrix} K(yf) - Y(yf, x\wedge2) \\ L(yf) - Y(yf, x*y) \\ M(yf) - Y(yf, x*z) \\ N(yf) - Y(yf, y\wedge2) \\ O(yf) - Y(yf, y*z) \\ P(yf) - Y(yf, z\wedge2) \end{bmatrix}$$

As previous, we can uniquely determine the unknown cubic coefficients present in the above matrix equation.

By comparing the terms of **degree = 2**, i.e. by collecting and equating the coefficients of x^2, x*y, x*z, y^2, y*z, z^2 in both expressions for **"yg"** we arrive at the six equations and we can collect these equations into a single matrix equation given below:

$$A^T \begin{bmatrix} q210 \\ 2*q120 \\ q111 \\ 3*q030 \\ 2*q021 \\ q012 \end{bmatrix} = \begin{bmatrix} K(yg) - Y(yg, x\wedge2) \\ L(yg) - Y(yg, x*y) \\ M(yg) - Y(yg, x*z) \\ N(yg) - Y(yg, y\wedge2) \\ O(yg) - Y(yg, y*z) \\ P(yg) - Y(yg, z\wedge2) \end{bmatrix}$$

By comparing the terms of **degree = 2**, i.e. by collecting and equating the coefficients of x^2, x*y, x*z, y^2, y*z, z^2 in both expressions for **"yh"** we arrive at the six equations and we can collect these equations into a single matrix equation given below:

$$A^T \begin{bmatrix} q201 \\ q111 \\ 2*q102 \\ q021 \\ 2*q012 \\ 3*q003 \end{bmatrix} = \begin{bmatrix} K(yh) - Y(yh, x\wedge 2) \\ L(yh) - Y(yh, x*y) \\ M(yh) - Y(yh, x*z) \\ N(yh) - Y(yh, y\wedge 2) \\ O(yh) - Y(yh, y*z) \\ P(yh) - Y(yh, z\wedge 2) \end{bmatrix}$$

As previous, we can uniquely determine the unknown cubic coefficients present in the above matrix equation.

By comparing the terms of **degree = 2**, i.e. by collecting and equating the coefficients of x^2, x*y, x*z, y^2, y*z, z^2 in both expressions for **"zf"** we arrive at the six equations and we can collect these equations into a single matrix equation given below:

$$A^T \begin{bmatrix} 3*r300 \\ 2*r210 \\ 2*r201 \\ r120 \\ r111 \\ r102 \end{bmatrix} = \begin{bmatrix} K(zf) - Y(zf, x\wedge 2) \\ L(zf) - Y(zf, x*y) \\ M(zf) - Y(zf, x*z) \\ N(zf) - Y(zf, y\wedge 2) \\ O(zf) - Y(zf, y*z) \\ P(zf) - Y(zf, z\wedge 2) \end{bmatrix}$$

As previous, we can uniquely determine the unknown cubic coefficients present in the above matrix equation.

By comparing the terms of **degree = 2**, i.e. by collecting and equating the coefficients of x^2, x*y, x*z, y^2, y*z, z^2 in both expressions for **"zg"** we arrive at the six equations and we can collect these equations into a single matrix equation given below:

$$A^T \begin{bmatrix} r210 \\ 2*r120 \\ r111 \\ 3*r030 \\ 2*r021 \\ r012 \end{bmatrix} = \begin{bmatrix} K(zg)-Y(zg,x\hat{}2) \\ L(zg)-Y(zg,x*y) \\ M(zg)-Y(zg,x*z) \\ N(zg)-Y(zg,y\hat{}2) \\ O(zg)-Y(zg,y*z) \\ P(zg)-Y(zg,z\hat{}2) \end{bmatrix}$$

As previous, we can uniquely determine the unknown cubic coefficients present in the above matrix equation.

By comparing the terms of **degree = 2**, i.e. by collecting and equating the coefficients of x^2, x*y, x*z, y^2, y*z, z^2 in both expressions for **"zh"** we arrive at the six equations and we can collect these equations into a single matrix equation given below:

$$A^T \begin{bmatrix} r201 \\ r111 \\ 2*r102 \\ r021 \\ 2*r012 \\ 3*r003 \end{bmatrix} = \begin{bmatrix} K(zh)-Y(zh,x\hat{}2) \\ L(zh)-Y(zh,x*y) \\ M(zh)-Y(zh,x*z) \\ N(zh)-Y(zh,y\hat{}2) \\ O(zh)-Y(zh,y*z) \\ P(zh)-Y(zh,z\hat{}2) \end{bmatrix}$$

As previous, we can uniquely determine the unknown cubic coefficients present in the above matrix equation.

(4) Instead of assuming BCW form for the chosen polynomials, $f, g, h$ if we will take initially the polynomials $f, g, h$ in which quadratic terms are present and as was done in (1), (2), (3) above for the cases of **degree = 0, 1, 2** if we will build the equations for **degree = 3**, i.e. by collecting and equating the coefficients of x^3, x^2*y, x^2*z, x*y^2, x*y*z, x*z^2, y^3, y^2*z, y*z^2, z^3 in both expressions for **“xf”** then we arrive at **ten equations** and we can collect these **ten equations** into a single matrix equation.

$$C^T \begin{bmatrix} 4*p400 \\ 3*p310 \\ 3*p301 \\ 2*p220 \\ 2*p211 \\ 2*p202 \\ p130 \\ p121 \\ p112 \\ p103 \end{bmatrix} = \begin{bmatrix} K_1(xf) - Y(xf, x\text{^}3) \\ K_2(xf) - Y(xf, x\text{^}2*y) \\ K_3(xf) - Y(xf, x\text{^}2*z) \\ K_4(xf) - Y(xf, x*y\text{^}2) \\ K_5(xf) - Y(xf, x*y*z) \\ K_6(xf) - Y(xf, x*z\text{^}2) \\ K_7(xf) - Y(xf, y\text{^}3) \\ K_8(xf) - Y(xf, y\text{^}2*z) \\ K_9(xf) - Y(xf, y*z\text{^}2) \\ K_{10}(xf) - Y(xf, z\text{^}3) \end{bmatrix}$$

We can easily check that **the components of the vector on right hand side of the above matrix equation are all equal to zero**. For example,

$$K_1(xf) = (b110*c201 + b210*c101) - (b101*c210 + b201*c110)$$

The $K_1(xf)$ above is obtained by collecting the coefficients of x^3 from xf=(gy*hz-gz*hy) and one observes that each term there contains some coefficient of quadratic term in $f$ or $g$ or $h$ which is already equal to zero due to chosen BCW form for the cubic polynomials $f, g, h$. Also,

$$\begin{aligned} Y(xf, x\text{^}3) &= 2*p200*a300 + p110*a300 + p101*c300 \\ &+ 3p300*2*a100*a200 + 2*p210(a100*b200 + a200*b100) \\ &+ 2*p201*(a100*c200 + a200*c100) + p120*2*b100*b200 \\ &+ p111*(b100*c200 + b200*c100) + p102*2*c100*c200 \end{aligned}$$

The $Y(xf, x\text{^}3)$ above is obtained by collecting the coefficients of x^3 from xf given in equation (5.6.24). Here also each term either contains some coefficient of quadratic term in $f$ or $g$ or $h$ which is already equal to zero due to chosen BCW form for the cubic polynomials $f, g, h$, or it contains some coefficient of quadratic term in $x$ or $y$ or $z$ as power series in $f, g, h$ which are also already equal to zero due to $p_{ijk} = 0, q_{ijk} = 0, r_{ijk} = 0$ whenever $i + j + k = 2$ as already seen on page 60 above when given polynomials have BCW form. It is easy to check that the ten by ten matrix $C^T$ is invertible, $\det(C^T) = 1$, and let $(C^T)^{-1} = D^T$, where $D^T$ can be obtained either by inverting $C^T$ or by forming inverse relations (again ten equations) and collecting them into a single matrix equation.

This procedure seen above for **"xf"** can be carried out for all other partial derivatives **"xg, xh, yf, yg, yh, zf, zg, zh"** and each time we will get a set of ten equations as above and each time we can collect those corresponding **ten equations** into a single matrix equation as above and we can show by same way that **the components of vector on right hand side of each such matrix equation are all equal to zero**. Thus, we can show that $p_{ijk} = 0, q_{ijk} = 0, r_{ijk} = 0$ whenever $i + j + k = 4$. By continuing with the same procedure for **degree = 4, 5, 6,… cases** we can further show that actually, $p_{ijk} = 0, q_{ijk} = 0, r_{ijk} = 0$ whenever $i + j + k \geq 4$.
Thus, three variable case is complete and now. We have seen that we can easily construct the inverse polynomials $x, y, z$ in terms of $f, g, h$ as basic variables and all the coefficients of terms in the inverse polynomials $x, y, z$ are in the terms of the coefficients of given polynomials $f, g, h$. Thus we have shown for the three variable case that the inverse functions exist and they are actually polynomials and we can construct them when the Jacobian is a nonzero constant (=1, say).

**6. The Several Variables Case:** We will be very brief in this section. We will just state the results that one can obtain by proceeding on similar lines for two and three variables cases.

(1) We can start with polynomials $u_i, i = 1, \cdots, n$, each one in the variables $x_i, i = 1, \cdots, n$ satisfying the Jacobi condition, namely,

$$J_x(u) = \det\left(\frac{\partial u_i}{\partial x_j}\right) = M = 1 \neq 0 \qquad (6.1)$$

We can apply inverse function theorem and obtain the matrix equation

$$\begin{pmatrix} \frac{\partial x_1}{\partial u_1} & \frac{\partial x_1}{\partial u_2} & \cdots & \frac{\partial x_1}{\partial u_n} \\ \frac{\partial x_2}{\partial u_1} & \frac{\partial x_2}{\partial u_2} & \cdots & \frac{\partial x_2}{\partial u_n} \\ \vdots & \vdots & & \vdots \\ \frac{\partial x_n}{\partial u_1} & \frac{\partial x_n}{\partial u_2} & \cdots & \frac{\partial x_n}{\partial u_n} \end{pmatrix} = \begin{pmatrix} \frac{\partial u_1}{\partial x_1} & \frac{\partial u_1}{\partial x_2} & \cdots & \frac{\partial u_1}{\partial x_n} \\ \frac{\partial u_2}{\partial x_1} & \frac{\partial u_2}{\partial x_2} & \cdots & \frac{\partial u_2}{\partial x_n} \\ \vdots & \vdots & & \vdots \\ \frac{\partial u_n}{\partial x_1} & \frac{\partial u_n}{\partial x_2} & \cdots & \frac{\partial u_n}{\partial x_n} \end{pmatrix}^{-1} = (W)^{-1} \qquad (6.2)$$

(2) Thus, we start formally with the given n polynomials, each in n variables, like,

$$u_i = \sum_{m_1=0}^{l_1^i} \sum_{m_2=0}^{l_2^i} \cdots \sum_{m_n=0}^{l_n^i} a^i{}_{m_1 m_2 \cdots m_n} x_1^{m_1} x_2^{m_2} \cdots x_n^{m_n} \qquad (6.3)$$

and formally construct the inverse functions (power series) like

$$x_i = \sum_{m_1=0} \sum_{m_2=0} \cdots \sum_{m_n=0} b^i{}_{m_1 m_2 \cdots m_n} u_1^{m_1} u_2^{m_2} \cdots u_n^{m_n} \qquad (6.4)$$

(3) From the matrix equation (6.2) we construct $n^2$ equations by taking formal partial derivative $\left(\frac{\partial x_i}{\partial u_j}\right)$ by using equation (6.4) involving terms containing unknown and to be determined coefficients and equating these partial derivatives as shown below using given polynomials defined in equations (6.3), viz,

$$\left(\frac{\partial x_i}{\partial u_j}\right) = (-1)^{(i+j)} \det(U_{ji}) \qquad (6.5)$$

where $\det(U_{ji}) =$ Determinant of the matrix obtained by deleting $j$-th row and $i$-th column of matrix $W$ whose inverse is given in the equation (6.2) above.

(4) By proceeding exactly as previous two and three variable cases, as they are applicable word to word; we determine the inverse functions and show easily that they are actually polynomials!!

(5) Considering the linear polynomials $u_i^*$, as is done in section 2.3, we can check the validity of the result like theorem 2.3.1, and thus can establish the invertible nature of matrices **like** $U_n^T$ for several variables case.

(6) We can obtain on similar lines the principle as well as the derived Jacobi conditions of **all types** for several variables.
The principle Jacobi condition here becomes

$$J_{(x)}(u) = \det\begin{pmatrix} a^1_{10\cdots0} & \cdots & a^1_{00\cdots1} \\ \vdots & & \vdots \\ a^n_{10\cdots0} & \cdots & a^n_{00\cdots1} \end{pmatrix} = M = 1 \neq 0$$

(7) As is done in section 4 we can establish the nonsingular nature of matrices **like** $U_n^T$ for several variables case, and can develop a theorem like theorem 4.1, namely, the determinant of the matrix **like** $U_n^T$, $n = 2,3,\cdots$ for several variables case is power of the corresponding Jacobian $M$, i.e. $\det(U_n^T) = M^k$ for some positive integer $k$. Moreover, the values $k$ considered in succession form a line ((**n+1)-th diagonal line**) in the Pascal triangle.

We have thus arrived at a position now to state the important result that we have obtained in this paper. The important result due to H. Bass, E. H. Connell, and D. Wright in [3] very much reduces the computational burden by their **important reduction** in degree achieved for the Jacobian problem. According to their result it is enough to settle the Jacobian problem for the special homogeneous polynomials $u \equiv u(x)$ of the special cubic form, namely, $u \equiv u(x) = x - H(x)$, where $H(tx) = t^3 H(x)$ for all $t \in k$ and all $x \in k^n$, $k$ being the ground field of characteristic zero. The procedures discussed in the paper can be easily applied to the cubic polynomials of this special form and one can show the vanishing of the "residuum" for these special type of polynomials as required, with much ease. Their result essentially is as follows:

**Theorem [BCW]**: The Jacobian conjecture is true for polynomials $u(x)$ having every number of variables $n$, and for every degree if and only if it is true for polynomials having every number of variables $n$, and having cubic degree i.e. having special kind of cubic-homogeneous form: $u(x) = x - H(x)$, where $H(\alpha x) = \alpha^3 H(x)$ for every $\alpha \in k$, the ground field of characteristic zero.

What we have essentially achieved in this paper is the following result:

**Theorem [Jacobian Conjecture]:** Given $n$ polynomials, each one of having special kind of cubic-homogeneous form: $u(x) = x - H(x)$, where $H(\alpha x) = \alpha^3 H(x)$ for every $\alpha \in k$, the ground field of characteristic zero. Thus, $u = (u_1, u_2, \cdots, u_n)$ are polynomials of special type mentioned above in $n$ variables $x = (x_1, x_2, \cdots, x_n)$ and their Jacobian $J_x(u) = \det\left(\frac{\partial u_i}{\partial x_j}\right)$, $i, j = 1,2,\cdots,n$ is a nonzero constant (=1) in the ground field $k$ of characteristic zero then we get a polynomial inverse for this system of polynomials i.e. we get each $x_i$ as special kind of cubic-homogeneous polynomial: $x(u) = u - H(u)$, where $H(\alpha u) = \alpha^3 H(u)$ for every $\alpha \in k$, the ground field of characteristic zero polynomial in $n$ variables $u = (u_1, u_2, \cdots, u_n)$.

□

# References

1. S. S. Abhyankar, Expansion Techniques in Algebraic Geometry, Tata Institute of Fundamental Research, Mumbai, 1977.
2. O. H. Keller, Ganze Cremona -Transformationen, Monatsh. Math. Physik, 47, 299-306, 1939.
3. H. Bass, E. H. Connell, and D. Wright, The Jacobian Conjecture: Reduction of Degree and Formal Expansion of the Inverse. Bull. Amer. Math. Soc., 7, 287-330, 1982.